\documentclass[11pt]{article}
\pdfoutput=1
\usepackage[utf8]{inputenc}
\usepackage{amsmath,amssymb,epsfig,bbm}
\usepackage{comment}
\usepackage[dvips]{color}
\usepackage[T1]{fontenc}
\newtheorem{defi}{Definition}
\newtheorem{prop}[defi]{Proposition}
\newtheorem{theo}[defi]{Theorem}
\newtheorem{conj}[defi]{Conjecture}
\newtheorem{lemm}[defi]{Lemma}
\newtheorem{coro}[defi]{Corollary}
\newtheorem{rema}[defi]{Remark}
\newtheorem{exem}[defi]{Example}
\newtheorem{exems}[defi]{Examples}

\newcommand{\bdefi}{\begin{defi}}
\newcommand{\edefi}{\end{defi}}
\newcommand{\bprop}{\begin{prop}}
\newcommand{\eprop}{\end{prop}}
\newcommand{\btheo}{\begin{theo}}
\newcommand{\etheo}{\end{theo}}
\newcommand{\blemm}{\begin{lemm}}
\newcommand{\brema}{\begin{rema}}
\newcommand{\erema}{\end{rema}}
\newcommand{\bexer}{\begin{exem}}
\newcommand{\eexer}{\end{exem}}
\newcommand{\bexems}{\begin{exems}}
\newcommand{\eexems}{\end{exems}}
\newcommand{\bconj}{\begin{conj}}
\newcommand{\econj}{\end{conj}}
\newcommand{\elemm}{\end{lemm}}
\newcommand{\bcoro}{\begin{coro}}
\newcommand{\ecoro}{\end{coro}}
\newcommand{\dem}{\noindent{\bf Proof. }}

\usepackage{mathrsfs}
\renewcommand\mathcal{\mathscr}

\newcommand{\M}{{\cal M}}
\newcommand{\N}{{\cal N}}

\newcommand{\V}{{\cal V}}
\newcommand{\W}{{\cal W}}

\renewcommand{\H}{{\cal H}}
\newcommand{\OOO}{{\cal O}}
\newcommand{\C}{{\cal C}}
\newcommand{\I}{{\cal I}}

\renewcommand{\P}{{\cal P}}

\newcommand{\maths}[1]{{\mathbb #1}}  

\newcommand{\RR}{\maths{R}}
\newcommand{\NN}{\maths{N}}
\newcommand{\CC}{\maths{C}}
\newcommand{\QQ}{\maths{Q}}

\newcommand{\SSS}{\maths{S}}

\newcommand{\HH}{\maths{H}}

\newcommand{\ZZ}{\maths{Z}}
\newcommand{\PP}{\maths{P}}

\newcommand{\aaa}{{\mathfrak a}}

\newcommand{\ppp}{{\mathfrak p}}

\newcommand{\mmm}{{\mathfrak m}}

\newcommand{\ra}{\rightarrow}
\newcommand{\bs}{\backslash}

\newcommand{\ov}[1]{{\overline #1}} 
\newcommand{\wt}[1]{{\widetilde{#1}}}

\newcommand{\ga}{\gamma}
\newcommand{\Ga}{\Gamma}

\newcommand{\cqfd}{\hfill$\Box$}

\newcommand{\card}{{\operatorname{Card}}}
\renewcommand{\Re}{{\operatorname{Re}}}
\renewcommand{\Im}{{\operatorname{Im}}}
\newcommand{\Vol}{\operatorname{Vol}}

\newcommand{\covol}{\operatorname{Covol}}
\newcommand{\id}{\operatorname{id}}

\newcommand{\PSL}{\operatorname{PSL}}

\newcommand{\dvol}{d\operatorname{Vol}}
\newcommand{\Det}{\operatorname{Det}}
\newcommand{\CAT}{\operatorname{CAT}}
\newcommand{\cinf}{\operatorname{C}^\infty}
\newcommand{\Perp}{\operatorname{Perp}}
\newcommand{\curv}{\operatorname{curv}_{\rm S}}

\newcommand{\hdr}{{\HH}^2_\RR}
\newcommand{\htr}{{\HH}^3_\RR}
\newcommand{\hcr}{{\HH}^5_\RR}
\newcommand{\hkr}{{\HH}^k_\RR}
\newcommand{\hnr}{{\HH}^n_\RR}

\newcommand{\SLO}{\operatorname{SL}_{2}(\OOO)}
\newcommand{\SLOK}{\operatorname{SL}_{2}(\OOO_K)}
\newcommand{\SLZ}{\operatorname{SL}_{2}(\ZZ)}

\newcommand{\PSLOK}{\operatorname{PSL}_{2}(\OOO_K)}
\newcommand{\PSLZ}{{\operatorname{PSL}_{2}(\ZZ)}}
\newcommand{\SLH}{\operatorname{SL}_{2}(\HH)}
\newcommand{\PSLH}{\operatorname{PSL}_{2}(\HH)}

\newcommand{\const}{\iota}
\newcommand{\discr}{\Delta}

\newcommand{\tr}{\operatorname{tr}}
\newcommand{\n}{\operatorname{n}}

\newcommand{\flow}[1]{{g^{#1}}}  
  
\newcommand{\mus}[1]{{\mu^{\rm s}_{#1}}}

\newcommand{\normal}[1]{\partial^1_{+}{#1}}

\newcounter{fig}

\title{Counting arcs in negative curvature}
\author{Jouni Parkkonen \and Fr\'ed\'eric Paulin} 
\date{}
\begin{document}
\bibliographystyle{../alphanum}
\maketitle

\begin{abstract} Let $M$ be a complete Riemannian manifold with
  negative curvature, and let $C_-,C_+$ be two properly immersed
  closed convex subsets of $M$. We survey the asymptotic behaviour of
  the number of common perpendiculars of length at most $s$ from $C_-$
  to $C_+$, giving error terms and counting with weights, starting
  from the work of Huber, Herrmann, Margulis and ending with the works
  of the authors. We describe the relationship with counting problems
  in circle packings of Kontorovich, Oh, Shah. We survey the tools
  used to obtain the precise asymptotics (Bowen-Margulis and Gibbs
  measures, skinning measures). We describe several arithmetic
  applications, in particular the ones by the authors on the
  asymptotics of the number of representations of integers by binary
  quadratic, Hermitian or Hamiltonian forms. 
\end{abstract}

\section{Introduction}
\label{sec:intro}

Let $M$ be a complete connected Riemannian manifold with negative
sectional curvature. Let $C_-$ and $C_+$ be two properly immersed
closed convex subsets of $M$ (see the end of \S
\ref{subsec:convexsubset} for a precise definition). For instance,
$C_-$ and $C_+$ could be points, totally geodesic immersed
submanifolds, Margulis cusp neighbourhoods, or images in $M$ of convex
hulls in a universal Riemannian cover of $M$ of limit sets of
subgroups of the fundamental group of $M$.  A {\em common
  perpendicular} from $C_-$ to $C_+$ is a locally geodesic path
$c:[a,b]\ra M$ such that $\dot c(a)$ is an outer unit normal vector to
$C_-$ and $-\dot c(b)$ is an outer unit normal vector to $C_+$ (see \S
\ref{subsec:convexsubset} and \S \ref{subsec:countingproblem} for
precise definitions when the boundary of $C_\pm$ is not smooth). The
(appropriately indexed) set of lengths of these common perpendiculars
is called the {\it (marked) ortholength spectrum} of $(C_-,C_+)$, and
variations on it have been introduced in particular cases of constant
curvatures by Basmajian, Bridgeman, Bridgeman-Kahn,
Martin-McKee-Wambach, Meyerhoff, Mirzakhani (see \S
\ref{subsec:countingproblem} for precisions).

The aim of this survey is to present several results on the asymptotic
behaviour as $s$ tends to $+\infty$ of the number $\N(s)$ of common
perpendiculars (counted with multiplicities, see \S
\ref{subsec:countingproblem}) between $C_-$ and $C_+$ with length at
most $s$.

We describe the first results on this problem by Huber \cite{Huber59},
Herrmann \cite{Herrmann62} and Margulis \cite{Margulis69}, see also
the surveys \cite{Babillot02a, Sharp04} and their references.  We
explain in \S \ref{sec:counting} how several works of
Duke-Rudnick-Sarnak on counting integral points on hyperboloids and of
Kontorovich, Oh and Shah on counting problems in circle and sphere
packings are related to counting problems of common perpendiculars. We
particularly emphasize the arithmetic applications developped in
\cite{ParPau12JMD,ParPau11BLMS,ParPau12ANT} to the counting of the
representations of integers by quadratic, Hermitian or quaternionic
binary forms, see \S \ref{sec:arithmapplic}. 

A few results are new, in particular the deduction of the
equidistribution of the outer unit normal bundles of equidistant
hypersurfaces of totally geodesic submanifolds in a hyperbolic
manifold from Eskin-McMullen's \cite{EskMcMul93} work (see \S
\ref{sec:eskinmcmullen}), and the computations of the constant relating
the Bowen-Margulis measure and the Liouville measure in constant
curvature and finite volume, as well as the one relating the skinning
measure and the Riemannian measure on the unit normal bundle of a
totally geodesic submanifold in constant curvature and finite volume
(see \S \ref{sec:finitevolhypman}).

We survey the main tools used for the counting results: the geometry
of negatively curved manifolds in \S \ref{sec:geometry} on one end,
and on the other hand, in \S \ref{sec:skinning}, the various measures,
as the Patterson-Sullivan densities, the Bowen-Margulis measures, the
skinning measures, that are needed to explicit the multiplicative
constant in front of the exponential term in the asymptotics of the
number $\N(s)$ of common perpendiculars. Skinning measures have first
been introduced by Oh-Shah in constant curvature for immersed balls,
horoballs and totally geodesic submanifolds, and are developped in general
in \cite{ParPau13a}.

We give in \S \ref{sec:mainstatement} a sketch of the proof of
the main counting result of \cite{ParPau13b}, which seems to contain
as particular cases all previous results on the asymptotics of
$\N(s)$, and to give many new ones. We conclude the paper by studying the
error term to this asymptotic equivalent under hypotheses of
exponential mixing of the geodesic flow (see \S
\ref{sec:errorterm}), and by giving counting results when weights
have been added to the common perpendiculars, by means of a
potential function, the main tools being then the Gibbs measures (see
\S \ref{sec:gibbs}).

To keep this paper to a reasonable length, we have chosen not to
develop the related counting problems of closed geodesics with lengths
at most $s$, which have been studied extensively (see for instance
\cite{Bowen72b,ParPol83,Parry84,Roblin03,PauPolSha}, and the surveys
\cite{Babillot02a,Sharp04} and their references).

\medskip {\small {\em Acknowledgment. } The authors thank the FIM of
  the ETH in Z\"urich, where most of this paper has been written, in
  particular for the support of the first author during the year
  2011-2012. We thank Hee Oh for the instigation to write \S
  \ref{sec:eskinmcmullen}.}

\section{Geometry and dynamics in negative 
curvature}    
\label{sec:geometry}

In this section, we survey briefly the required background on the
geometry and dynamics of negatively curved Riemannian manifolds,
considered as locally $\CAT(-\kappa)$ spaces (using for instance
\cite{BriHae99} as a general reference), with a particular emphasis on
the metric aspects and the regularity properties.

For every $n\geq 2$, we denote by $\HH^n_\RR$ %any model of 
the real
hyperbolic space of dimension $n$ and constant sectional curvature
$-1$.

For every $\epsilon>0$, we denote by $\N_\epsilon A$ the closed
$\epsilon$-neighbourhood of a subset $A$ of a metric space, and by
convention $\N_0A=\overline{A}$. Recall that a map $f:X\ra Y$ between
two metric spaces is called $\alpha$-{\em H\"older}-continuous, where
$\alpha\in\,]0,1]$, if there exists $c,c'>0$ such that for every $x,y$
in $X$ with $d(x,y)\leq c'$, we have $d(f(x),f(y))\leq
c\,d(x,y)^\alpha$, and {\em H\"older}-continuous if there exists
$\alpha\in\,]0,1]$ such that $f$ is $\alpha$-H\"older-continuous.

\subsection{Geometry of the unit tangent bundle}
\label{subsec:geomunittgbun}

We denote by $\pi: TN\to N$ the tangent bundle of any (smooth)
Riemannian manifold $N$, and again by $\pi: T^1 N\to N$ its unit
tangent bundle. Recall that the Levi-Civita connexion $\nabla$ of $N$
gives a decomposition $TTN=V\oplus H$ of the vector bundle $TTN\ra TN$
into the direct sum of two smooth vector subbundles $V\ra TN$ and
$H\ra TN$, called vertical  and horizontal, such that if $\pi_V:TTN\ra
V$ is the linear projection of $TTN$ onto $V$ parallelly to $H$, if
$H_v$ and $V_v$ are the fibers of $H$ and $V$ above $v\in TN$, then

$\bullet$~~we have $V_v=\operatorname{Ker} T_v\pi=T_v(T_{\pi(v)}N)=
T_{\pi(v)}N$;

$\bullet$~~the restriction $T\pi_{\mid H_v}:H_v\ra T_{\pi(v)}N$
of the tangent map of $\pi$  to $H_v$ is a linear isomorphism;

$\bullet$~~for every smooth vector field $X:N\ra TN$ on $N$, we have
$\nabla_v X=\pi_V\circ TX(v)$.

The manifold $TN$ has a unique Riemannian metric, called {\em Sasaki's
  metric}, such that for every $v\in TM$, the map $T\pi_{\mid H_v} :
H_v\ra T_{\pi(v)}N$ is isometric, the restriction to $V_v$ of Sasaki's
scalar product is the Riemannian scalar product on $T_{\pi(v)}N$, and
the decomposition $T_vTN=V_v\oplus H_v$ is orthogonal. We endow the
smooth submanifold $T^1 N$ of $TN$ with the induced Riemannian metric,
also called {\em Sasaki's metric}.  The fiber $T^1_xN$ of every $x\in
N$ is then isometric to the standard unit sphere $\SSS^{n-1}$ of the
standard Euclidean space $\RR^n$, if $n$ is the dimension of $N$.

The Riemannian measure $\dvol_{T^1N}$ of $T^1N$, called {\em
  Liouville's measure}, disintegrates under the fibration $\pi:T^1N
\ra N$ over the Riemannian measure $\dvol_{N}$ of $N$, as
$$
\dvol_{T^1N}=\int_{x\in N} \;\dvol_{T^1_xN}\;\dvol_{N}(x)\;,
$$
where $\dvol_{T^1_xN}$ is the spherical measure on the fiber $T^1_xN$
of $\pi$ above $x\in N$. In particular,
$$
\Vol(T^1N)=\Vol(\SSS^{n-1})\Vol(N)\;.
$$

\subsection{H\"older structure on the boundary at
  infinity}
\label{subsec:holderstrcut}

Let $\wt M$ be a complete simply connected Riemannian manifold with
dimension at least $2$ and sectional curvature at most $-1$, and let
$x_0\in\wt M$. To shorten the exposition, we assume in this survey
that $\wt M$ has pinched negative sectional curvature $-b^2\leq K\leq
-1$ (where $b\in\;[1,+\infty[$), though this is not necessary except
when working with Gibbs measures in \S \ref{sec:gibbs}, see
\cite{ParPau13a,ParPau13b} for the extensions. The error term
estimates of \S \ref{sec:errorterm} require another geometric
assumption, that the sectional curvature of $M$ has bounded
derivatives.

We denote by $\partial_{\infty}\wt M$ the boundary at infinity of $\wt
M$, with its usual H\"older structure when the sectional curvature of
$M$ has bounded derivatives and its usual conformal structure, which
we describe below.  Recall that a {\em H\"older structure} on a
topological manifold $X$ is a maximal atlas of charts $(U,\varphi)$,
where $\varphi: U\ra \varphi(U)$ is a homeomorphism between an open
subset $U$ of $X$ and an open subset of a fixed smooth manifold, such
that the transition maps are $\alpha$-H\"older homeomorphisms for some
$\alpha>0$.

Two geodesic rays $\rho,\rho':[0,+\infty[\;\ra\wt M$ are {\em
  asymptotic} if their images are at finite Hausdorff distance, or
equivalently if there exists $c>0$ and $t_0\in\RR$ such that the
inequality $d(\rho(t),\rho'(t+t_0))\leq c\;e^{-t}$ holds for all $t\in
[\max\{0,-t_0\}, +\infty[\,$. The {\em boundary at infinity}
$\partial_{\infty}\wt M$ of $\wt M$ is the quotient topological space
of the space of geodesic rays, endowed with the compact-open topology,
by the equivalence relation ``to be asymptotic to''.  The asymptotic
class of a geodesic ray is called its {\em point at infinity}. For all
$x\in M$ and $\xi\in \partial_{\infty} \wt M$, there exists a unique
geodesic ray with origin $x$ and point at infinity $\xi$, whose image
we denote by $[x,\xi[\,$.  Given two distinct points at infinity
$\xi,\eta \in \partial_{\infty} \wt M$, there exists a unique (up to
translation on the source) geodesic line $\rho:\RR\ra\wt M$ such that
the points at infinity of the geodesic rays $t\mapsto \rho(-t)$ and
$t\mapsto \rho(t)$, $t\in[0,+\infty[$, are $\xi$ and $\eta$,
respectively.  We denote the image of such a geodesic line by $]\xi,\eta[\,$.

For every $x\in \wt M$, the map $\theta_x$ from $T^1_x\wt M$ to
$\partial_{\infty}\wt M$, which sends $v\in T^1_x\wt M$ to the point
at infinity of the geodesic ray with tangent vector at the origin $v$,
is a homeomorphism. We define the angle $\angle_x(y,z)$ of two
geodesic segments or rays with the same origin $x$ and endpoints
$y,z\in (\wt M-\{x\})\cup \partial_{\infty}\wt M$ as the angle of
their tangent vectors at $x$. The disjoint union $\wt
M\cup \partial_{\infty}\wt M$ has a unique compact metrisable
topology, inducing the original topologies on $\wt M$ and on
$\partial_{\infty} \wt M$, such that a sequence of points $(y_n)_{n\in
  \NN}$ in $\wt M$ converges to a point $\xi\in \partial_{\infty} \wt
M$ if and only if $\lim_{n\ra+\infty} d(y_n,x_0)=+\infty$ and
$\lim_{n\ra+\infty} \angle_{x_0}(y_n,\xi) =0$. An isometry $\ga$ of
$\wt M$ uniquely extends to a homeomorphism of $\wt M
\cup \partial_{\infty} \wt M$, and we will also denote by $\ga$ its
extension to the boundary at infinity.

When the sectional curvature of $M$, besides being pinched negative,
has bounded derivatives, it is known since Anosov (see also
\cite{Brin95}, \cite[\S 7.1]{PauPolSha}) that the maps $\theta_x^{-1}
\circ \theta_{x'}: T^1_{x'} M \ra T^1_{x}M $, for all $x,x'\in \wt M$,
are $\alpha$-H\"older homeomorphisms for some $\alpha>0$. Hence there
is then a unique H\"older structure on the topological manifold
$\partial_{\infty}\wt M$ such that $\theta_{x}$ is a H\"older
homeomorphism, for every $x\in \wt M$. Furthermore, the isometries of
$\wt M$ are then $\alpha$-H\"older homeomorphisms of $\partial_{\infty}\wt
M$ for some $\alpha>0$.

\subsection{Conformal structure on the boundary at
  infinity}
\label{subsec:confstruc}

Let us now define the natural conformal structure on
$\partial_{\infty}\wt M$. Recall that two distances $d$ and $\delta$
on a set $Z$ are called {\em conformally equivalent} if they induce
the same topology and if for every $z_0\in Z$, the limit $\lim_{x\ra
  z_0, x\neq  z_0} \frac{d(x,z_0)}{\delta(x,z_0)}$ exists and is
strictly positive.  The relation ``to be conformally equivalent to''
is an equivalence relation on the set of distances on $Z$, and a {\em
  conformal structure} on $Z$ is an equivalence class thereof.

Let $Z$ and $Z'$ be two sets endowed with a conformal structure, and
let $d$ and $d'$ be distances on $Z$ and $Z'$ representing them. A
bijection $\ga:Z\ra Z'$ is {\em conformal} if the distances $d$ and
$\ga^*d':(x,y)\mapsto d'(\ga x,\ga y)$ are conformally equivalent.
%, that is, if for every $z_0\in Z$, the limit $\lim_{x\ra z_0, x\neq
%  z_0} \frac{d'(\ga x,\ga z_0)}{d(x,z_0)}$ exists and is strictly positive. 
This does not depend on the choice of representatives $d$
and $d'$ of the conformal structures of $Z$ and $Z'$.

The {\em Busemann cocycle} of $\wt M$ is the continuous map
$\beta\colon \wt M\times\wt M\times\partial_{\infty} \wt M\to\RR$
defined by
$$
(x,y,\xi)\mapsto \beta_{\xi}(x,y)=
\lim_{t\to+\infty}d(\rho_t,x)-d(\rho_t,y)\;,
$$
where $\rho:t\mapsto \rho_t$ is any geodesic ray with point at
infinity $\xi$. The above limit exists and is independent of
$\rho$. The Busemann cocycle is H\"older-continuous when the sectional
curvature of $M$ has bounded derivatives, and satisfies the
following equivariance and cocycle properties:
\begin{equation}\label{eq:cocycle}
\beta_{\ga \xi}(\ga x,\ga y)=\beta_{\xi}(x,y)\;\;\;{\rm and}\;\;\;
\beta_{\xi}(x,y)+\beta_{\xi}(y,z)=\beta_{\xi}(x,z)\;,
\end{equation}
for all $\xi\in\partial_{\infty}\wt M$, all $x,y,z\in\wt M$ and every
isometry $\ga$ of $\wt M$. In particular, $\beta_\xi(y,x)=
-\beta_{\xi}(x,y)$. By the triangular inequality, we have
\begin{equation}\label{eq:majobusema}
|\beta_{\xi}(x,y)|\leq d(x,y)\;.
\end{equation}
If $y$ is a point in the (image of the) geodesic ray from $x$ to
$\xi$, then $\beta_{\xi}(x,y)=d(x,y)$.  For every $y\in\wt M$ and
$\xi\in \partial_\infty \wt M$, the map $x\mapsto \beta_{\xi}(x,y)$ is
smooth and $1$-Lipschitz.

For every $\xi\in\partial_\infty \wt M$, the {\em horospheres centered
  at $\xi$} are the level sets of the map $y\mapsto \beta_\xi(y,x_0)$
from $\wt M$ to $\RR$, and the (closed) {\em horoballs centered at
  $\xi$} are its sublevel sets. Horoballs are closed (strictly) convex
subsets of $\wt M$. A horosphere centered at $\xi$ is a smooth
hypersurface of $\wt M$, orthogonal to the geodesic lines having $\xi$
as a point at infinity. In the upper halfspace model of the real
hyperbolic space $\HH^n_\RR$, the horospheres are the horizontal
affine hyperplanes therein or the Euclidean spheres therein tangent to
the horizontal coordinate hyperplane, with the point of tangency
removed. In the ball model of the real hyperbolic space $\HH^n_\RR$,
the horospheres (respectively horoballs) are the Euclidean spheres
(respectively balls) tangent to the unit sphere and contained in the
unit ball, with the point of tangency removed.

The horoballs are limits of big balls (and their centres the limits of
the centres thereof, explaining the terminology). More precisely, if
$\rho$ is a geodesic ray in $\wt M$ with point at infinity $\xi$, if
$B(t)$ is the ball of centre $\rho(t)$ and radius $t$, then the map
$t\mapsto B(t)$ converges to the horoball $H\!B$ of centre $\xi$
containing $\rho(0)$ in its boundary, for Chabauty's topology on
closed subsets of $\wt M$ (that is, for the Hausdorff convergence on
compact subsets: for every compact subset $K$ of $\wt M$, as $t\ra
+\infty$, the closed subset $(B(t)\cap K)\cup\,\overline{^cK}$
converges to the closed subset $(H\!B\cap K)\cup\,\overline{^cK}$ for
the Hausdorff distance).

\medskip
For every $x\in \wt M$, for all distinct $\xi,\eta\in\partial_\infty
\wt M$, the {\em visual distance} between $\xi$ and $\eta$ seen from
$x$ is
$$
d_x(\xi,\eta)= \lim_{t\ra+\infty} e^{-\frac{1}{2}\big(d(x,\,\rho_\xi(t))
  +d(x,\,\rho_\eta(t))-d(\rho_\xi(t),\,\rho_\eta(t))\big)}\;,
$$
for any geodesic rays $\rho_\xi$ and  $\rho_\eta$ with point at infinity
$\xi$ and $\eta$, respectively. Equivalently, with $t\mapsto \rho_t$
and $t\mapsto\rho'_t$ the geodesic rays with origin $x$ converging to
$\xi$ and $\eta$, we have
$$
d_x(\xi,\eta)= \lim_{t\ra+\infty} e^{\frac{1}{2}d(\rho_t,\,\rho'_t) -t}\;.
$$
Again equivalently, if $u$ is any point on the geodesic line between
$\xi$ and $\eta$, then
\begin{equation}\label{eq:defdistvistroi}
d_x(\xi,\eta)= e^{-\frac 12(\beta_\xi(x,\,u)+\beta_{\eta}(x,\,u))}\;.
\end{equation} 
Define $d_x(\xi,\eta)=0$ if $\xi=\eta$.

For every $x\in \wt M$, the above limits exist and the three formulas
coincide. The map $d_x:\partial_\infty \wt M\times \partial_\infty \wt
M\ra [0,+\infty[$ is a distance, inducing the original topology on
$\partial_\infty \wt M$. For all $x,y\in \wt M$, for all distinct
$\xi,\eta\in\partial_\infty \wt M$, and for every isometry $\ga$ of
$\wt M$, we have
$$
e^{-d(x,\;]\xi,\,\eta[)}\leq d_x(\xi,\eta)\leq 
(1+\sqrt{2})\,e^{-d(x,\;]\xi,\,\eta[)}\;,
$$
\begin{equation}\label{eq:rapportdistvisun}
\frac{d_x(\xi,\eta)}{d_y(\xi,\eta)}=
e^{-\frac 12(\beta_\xi(x,\,y)+\beta_\eta(x,\,y))}\;,
\end{equation}
\begin{equation}\label{eq:equivardistvis}
d_{\ga x}(\ga \xi,\ga \eta)=d_x(\xi,\eta)\;.
\end{equation}

\medskip

It follows from Equation \eqref{eq:rapportdistvisun} that the visual
distances $d_x$ for $x\in \wt M$ belong to the same conformal
structure on $\partial_\infty \wt M$. It follows from Equation
\eqref{eq:rapportdistvisun} and Equation \eqref{eq:equivardistvis}
that the (boundary extensions of) the isometries of $\wt M$ are
conformal bijections for this conformal structure. Furthermore, these
equations and Equation \eqref{eq:majobusema} imply that the isometries
of $\wt M$ are bilipschitz homeomorphisms for any visual distance: for
all $x\in \wt M$, for all distinct $\xi,\eta\in\partial_\infty \wt
M$, we have
$$
e^{-2\,d(x,\ga x)}\;d_x(\xi,\eta)\leq 
d_x(\ga\xi,\ga\eta)\leq 
e^{2\,d(x,\ga x)}\;d_x(\xi,\eta)\;.
$$

\subsection{Stable and unstable leaves of the geodesic flow}
\label{subsec:stabunstab}

We now turn to the description of the dynamics of the geodesic flow on
$\wt M$.

The unit tangent bundle $T^1N$ of a complete Riemannian manifold $N$
can be identified with the set of locally geodesic lines (parametrised
by arclength) $\ell\colon \RR\to N$ in $N$, endowed with the
compact-open topology. More precisely, we identify a locally geodesic
line $\ell$ and its (unit) tangent vector $\dot{\ell}(0)$ at time
$t=0$ and, conversely, any $v\in T^1 N$ is the tangent vector at time
$t=0$ of a unique locally geodesic line.  We will use this
identification without mention in this survey. In
particular, the base point projection $\pi\colon T^1 N\to N$ is given
by $\pi(\ell)=\ell(0)$.

The {\em geodesic flow} on $T^1N$ is the smooth $1$-parameter group
$(\flow{t})_{t\in\RR}$, where $\flow{t}\ell\,(s)=\ell(s+t)$, for all
$\ell\in T^1N$ and $s,t\in\RR$.  We denote by $\iota:T^1N\ra T^1N$ the
{\em antipodal (flip) map} $v\mapsto -v$. We have $\iota\circ \flow t=
\flow{-t}\circ \iota$. The isometry group of $N$ acts on the space of
geodesic lines in $N$ by postcomposition: $(\ga,\ell) \mapsto
\ga\circ\ell$, and this action commutes with the geodesic flow and the
antipodal map.

For every unit tangent vector $v\in T^1\wt M$, let $v_{-}=v(-\infty)$
and $v_{+}=v(+\infty)$ be the two endpoints in the sphere at infinity
of the geodesic line defined by $v$. Let $\partial_{\infty}^2\wt M$ be
the open subset of $\partial_{\infty}\wt M\times\partial_{\infty}\wt M
$ which consists of pairs of distinct points at infinity, with the
restriction of the product H\"older structure when the sectional
curvature of $M$ has bounded derivatives.  {\em Hopf's
  parametrisation} (see \cite{Hopf71}) of $T^1\wt M$ is the 
homeomorphism from $T^1\wt M$ to $\partial_{\infty}^2\wt M\times \RR$
sending $v\in T^1\wt M$ to the triple $(v_{-},v_{+},t)
\in \partial_{\infty}^2\wt M\times \RR$, where $t$ is the signed
(algebraic) distance of $\pi(v)$ from the closest point $p_{v,x_0}$ to
$x_{0}$ on the (oriented) geodesic line defined by $v$. Hopf's
  parametrisation is a H\"older homeomorphism when the sectional
curvature of $M$ has bounded derivatives. In this survey,
we will identify an element of $T^1\wt M$ with its image by Hopf's
parametrisation. The geodesic flow acts by $\flow s(v_{-},v_{+},t)
=(v_{-},v_{+},t+s)$ and, for every isometry $\ga$ of $\wt M$, the
image of $\ga v$ is $(\ga v_{-}, \ga v_{+}, t+ t_{\ga,v})$, where
$t_{\ga,v}$ is the signed distance from $\ga p_{v,x_0}$ to $p_{\ga
  v,x_0}$.  Furthermore, in these coordinates, the antipodal map
$\iota$ is $(v_-,v_+,t)\mapsto (v_+,v_-,-t)$.

\medskip
The {\em strong stable manifold} of $v\in T^1\wt M$ is 
$$
W^{\rm ss}(v)=\{v'\in T^1\wt M:d(v(t), v'(t))\to
0 \textrm{ as } t\to+\infty\},  
$$
and  the {\em strong unstable manifold} of $v$ is 
$$
W^{\rm su}(v)=\{v'\in T^1\wt M:d(v(t), v'(t))\to 0 
\textrm{ as } t\to-\infty\}, 
$$

\smallskip\noindent
\begin{minipage}{9.9cm} ~~~ The projections in $\wt M$ of the strong
  unstable and strong stable manifolds of $v\in T^1\wt M$, denoted by
  $H_{-}(v)=\pi(W^{\rm su}(v))$ and $H_{+}(v)=\pi(W^{\rm ss}(v))$, are
  called, respectively, the {\em unstable and stable horospheres} of
  $v$, and are the horospheres containing $\pi(v)$ centered at $v_{-}$
  and $v_{+}$, respectively. The unstable horosphere of $v$ coincides
  with the zero set of the map $x\mapsto f_{-}(x)= \beta_{v_-}(x,
  \pi(v))$, and, similarly, the stable horosphere of $v$ coincides
  with the zero set of $x\mapsto f_{+}(x)= \beta_{v_+}(x,\pi(v))$.
  The corresponding sublevel sets $H\!B_{-}(v)= f_{-}^{-1}(]-\infty,
  0])$ and $H\!B_{+}(v)= f_{+}^{-1}(]-\infty,0])$ are called the {\em
    unstable and stable horoballs} of $v$.
\end{minipage}
\begin{minipage}{5cm}
\begin{center}
\begin{picture}(0,0)%
\includegraphics{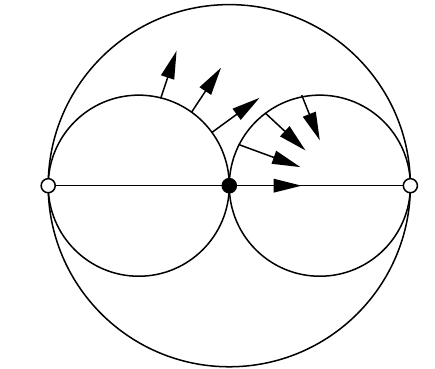}%
\end{picture}%
\setlength{\unitlength}{3812sp}%
\begingroup\makeatletter\ifx\SetFigFont\undefined%
\gdef\SetFigFont#1#2#3#4#5{%
  \reset@font\fontsize{#1}{#2pt}%
  \fontfamily{#3}\fontseries{#4}\fontshape{#5}%
  \selectfont}%
\fi\endgroup%
\begin{picture}(2145,1814)(1381,-1688)
\put(3511,-871){\makebox(0,0)[lb]{\smash{{\SetFigFont{11}{13.2}{\rmdefault}{\mddefault}{\updefault}{\color[rgb]{0,0,0}$v_+$}%
}}}}
\put(1396,-871){\makebox(0,0)[lb]{\smash{{\SetFigFont{11}{13.2}{\rmdefault}{\mddefault}{\updefault}{\color[rgb]{0,0,0}$v_-$}%
}}}}
\put(1891,-1366){\makebox(0,0)[lb]{\smash{{\SetFigFont{11}{13.2}{\rmdefault}{\mddefault}{\updefault}{\color[rgb]{0,0,0}$H_-(v)$}%
}}}}
\put(2746,-961){\makebox(0,0)[lb]{\smash{{\SetFigFont{11}{13.2}{\rmdefault}{\mddefault}{\updefault}{\color[rgb]{0,0,0}$v$}%
}}}}
\put(2746,-1366){\makebox(0,0)[lb]{\smash{{\SetFigFont{11}{13.2}{\rmdefault}{\mddefault}{\updefault}{\color[rgb]{0,0,0}$H_+(v)$}%
}}}}
\end{picture}%

\end{center}
\end{minipage}

\medskip The union for $t\in\RR$ of the images under $\flow t$ of the
strong stable manifold of $v\in T^{1}\wt M$ is the {\em stable
  manifold} $W^{\rm s}(v)=\bigcup_{t\in\RR}\flow t W^{\rm ss}(v)$ of
$v$, which consists of the elements $v'\in T^1\wt M$ with $v'_+=
v_+$. Similarly, the union of the images under the geodesic flow at
all times of the strong unstable manifold of $v$ is the {\em unstable
  manifold} $W^{\rm u}(v)$ of $v$, which consists of the elements
$v'\in T^1\wt M$ with $v'_-=v_-$.

The subspaces $W^{\rm ss}(v)$ and $W^{\rm su}(v)$ (which are the lifts
by the inner and outer unit normal vectors of the unstable and stable
horospheres of $v$, respectively), as well as $W^{\rm s}(v)$ and
$W^{\rm u}(v)$, are smooth submanifolds of $T^1\wt M$. The maps from
$\RR\times W^{\rm ss}(v)$ to $W^{\rm s}(v)$ defined by $(t,v')\mapsto
\flow tv'$ and from $\RR\times W^{\rm su}(v)$ to $W^{\rm u}(v)$
defined by $(t,v')\mapsto \flow tv'$ are smooth diffeomorphisms. We
have $W^{\rm ss}(\iota v)=\iota W^{\rm su}(v)$.

\bigskip
\noindent{\bf Hamenst\"adt's distance on stable and unstable leaves}

\medskip For every $v\in T^1\wt M$, let $d_{W^{\rm ss}(v)}$ be {\em
  Hamenst\"adt's distance} on the strong stable leaf of $v$, defined
as follows (see \cite{Hamenstadt89}, \cite[Appendix]{HerPau97}, as well
as \cite[\S 2.2]{HerPau10} for a generalisation when the horosphere
$H_{+}(v)$ is replaced by the boundary of any nonempty closed convex
subset): for all $w,w'\in W^{\rm ss}(v)$, we have
$$
d_{W^{\rm ss}(v)}(w,w') = 
\lim_{t\ra+\infty} e^{\frac{1}{2}d(w(-t),\;w'(-t))-t}\;.
$$
This limit exists, and Hamenst\"adt's distance is a distance inducing
the original topology on $W^{\rm ss}(v)$. For all $w,w'\in W^{\rm ss}
(v)$ and for every isometry $\ga$ of $\wt M$, we have
$$
d_{W^{\rm ss}  (\ga v)}(\ga w,\ga w')= d_{W^{\rm ss}(v)}(w,w')\;.
$$
For all $v\in T^1\wt M$, $s\in\RR$ and $w,w'\in W^{\rm ss} (v)$, we
have
$$
d_{W^{\rm ss} (\flow sv)}(\flow sw,\flow sw')=e^{-s}d_{W^{\rm ss}(v)}(w,w')\;.
$$

For every horosphere $H$ in $\wt M$ with center $H_\infty$, we also
have a distance $d_H$ on the open subset $\partial_\infty\wt
M-\{H_\infty\}$ defined by
$$
d_{H}(\xi,\eta) = \lim_{t\ra +\infty}e^t\;d_{\rho_t}(\xi,\eta)=
\lim_{t\ra+\infty} e^{\frac{1}{2}d(\xi_t,\;\eta_t)-t}\;,
$$
where $t\mapsto\rho_t$ is any geodesic ray with origin a point of $H$
and point at infinity $H_\infty$, and $t\mapsto\xi_t$ and $t\mapsto
\eta_t$ are the geodesic lines in $\wt M$ with origin $H_\infty$,
passing at time $t=0$ through $H$, and with endpoints $\xi$ and $\eta$,
respectively. Using the homeomorphism from $W^{\rm ss} (v)$ to
$\partial_\infty\wt M-\{v_+\}$ defined by $w\mapsto w_-$, we have
$$
d_{W^{\rm ss} (v)}(w,w')=d_{H_+(v)}(w_-,w'_-)\;.
$$
The distance $d_H$ and the restriction of any visual distance to
$\partial_\infty\wt M-\{H_\infty\}$ are conformally equivalent, since
for all $x\in H$ and $\xi,\eta\in\partial_\infty\wt M-\{H_\infty\}$,
with the above notation $t\mapsto\xi_t$ and $t\mapsto \eta_t$, we 
have
$$
\frac{d_{H}(\xi,\eta)}{d_x(\xi,\eta)}=
e^{-\frac 12(\beta_\xi(\xi_0,\,x)+\beta_\eta(\eta_0,\,x))}\;.
$$

\medskip
\noindent{\bf The Anosov property of the geodesic flow}

\medskip The strong stable manifolds, stable manifolds, strong
unstable manifolds and unstable manifolds are the (smooth) leaves of
continuous foliations on $T^1\wt M$, invariant under the geodesic flow
and the isometry group of $\wt M$, denoted by $\W^{\rm ss}, \W^{\rm
  s}, \W^{\rm su}$ and $\W^{\rm u}$, respectively. They are H\"older
foliations when the sectional curvature of $M$ has bounded
derivatives. When $\wt M$ is a symmetric space (that is, up to
homothety, when $\wt M$ is isometric to the real, complex,
quaternionic hyperbolic $n$-space or to the octonionic hyperbolic
plane), then the strong stable, stable, strong unstable and unstable
foliations are smooth. But in general, the H\"older regularity cannot
be much improved, as we will explain in \S \ref{sec:discrisomgroups}.

\begin{center} 
\begin{picture}(0,0)%
\includegraphics{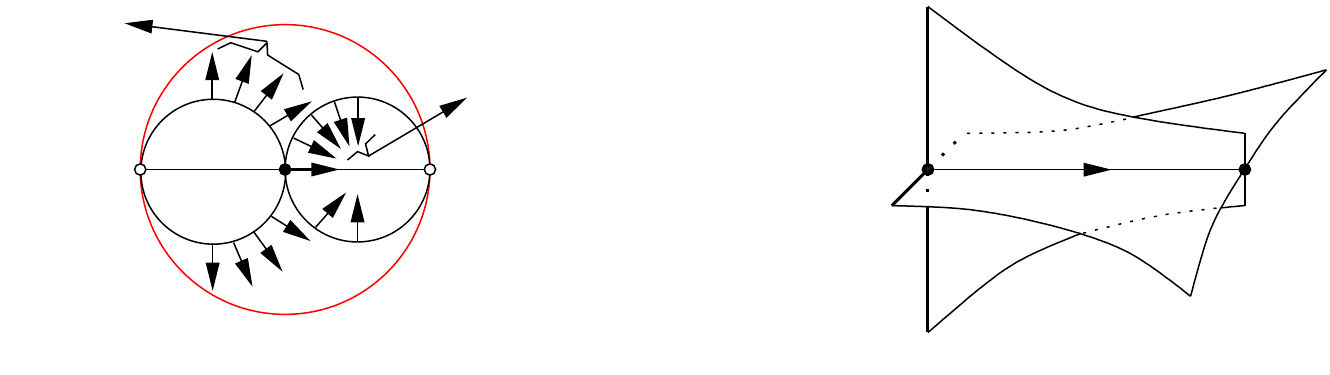}%
\end{picture}%
\setlength{\unitlength}{3812sp}%
\begingroup\makeatletter\ifx\SetFigFont\undefined%
\gdef\SetFigFont#1#2#3#4#5{%
  \reset@font\fontsize{#1}{#2pt}%
  \fontfamily{#3}\fontseries{#4}\fontshape{#5}%
  \selectfont}%
\fi\endgroup%
\begin{picture}(6604,1853)(3264,-1442)
\put(4754,-562){\makebox(0,0)[lb]{\smash{{\SetFigFont{11}{13.2}{\rmdefault}{\mddefault}{\updefault}{\color[rgb]{0,0,0}$v$}%
}}}}
\put(7726,-372){\makebox(0,0)[lb]{\smash{{\SetFigFont{11}{13.2}{\rmdefault}{\mddefault}{\updefault}{\color[rgb]{0,0,0}$v$}%
}}}}
\put(5611,-94){\makebox(0,0)[lb]{\smash{{\SetFigFont{11}{13.2}{\rmdefault}{\mddefault}{\updefault}{\color[rgb]{0,0,0}$W^{\rm ss}(v)$}%
}}}}
\put(7242,156){\makebox(0,0)[lb]{\smash{{\SetFigFont{11}{13.2}{\rmdefault}{\mddefault}{\updefault}{\color[rgb]{0,0,0}$W^{\rm ss}(v)$}%
}}}}
\put(7066,-601){\makebox(0,0)[lb]{\smash{{\SetFigFont{11}{13.2}{\rmdefault}{\mddefault}{\updefault}{\color[rgb]{0,0,0}$W^{\rm su}(v)$}%
}}}}
\put(3279,252){\makebox(0,0)[lb]{\smash{{\SetFigFont{11}{13.2}{\rmdefault}{\mddefault}{\updefault}{\color[rgb]{0,0,0}$W^{\rm su}(v)$}%
}}}}
\put(5239,-1157){\makebox(0,0)[lb]{\smash{{\SetFigFont{11}{13.2}{\rmdefault}{\mddefault}{\updefault}{\color[rgb]{0,0,0}$\wt M$}%
}}}}
\put(8717,-1378){\makebox(0,0)[lb]{\smash{{\SetFigFont{11}{13.2}{\rmdefault}{\mddefault}{\updefault}{\color[rgb]{0,0,0}$T^1\wt M$}%
}}}}
\put(9496,-511){\makebox(0,0)[lb]{\smash{{\SetFigFont{11}{13.2}{\rmdefault}{\mddefault}{\updefault}{\color[rgb]{0,0,0}$g^tv$}%
}}}}
\end{picture}%

\end{center}

Let $N=T^1\wt M$. The vector field $Z:N\ra TN$ defined by $v\mapsto
Z(v)=\frac{d}{dt} \flow{t}v$ is called the {\em geodesic vector
  field}. The geodesic flow $(\flow{t})_{t\in\RR}$ on the Riemannian
manifold $N$ is a {\em contact Anosov flow}. That is, the vector
bundle $TN\ra N$ is the direct sum of three topological vector
subbundles $TN=E_{\rm su}\oplus E_{0}\oplus E_{\rm ss}$ that are
invariant under $(\flow{t})_{t\in\RR}$, where $E_{0}\cap T_vN= \RR
Z(v)$, $E_{\rm su}\cap T_vN= T_vW^{\rm su}(v)$, $E_{\rm ss}\cap T_vN=
T_vW^{\rm ss}(v)$, and there exist two constants $c,\lambda>0$ such that
for every $t>0$, we have (see the above picture on the right)
$$
\|T_v\,\flow{t}\,_{\mid E_{\rm ss}}\|\leq c\,e^{-\lambda t}\;\;{\rm and}\;\;
\|T_v\,\flow{-t}\,_{\mid E_{\rm su}}\|\leq c\,e^{-\lambda t}\;.
$$
Furthermore, if $\alpha$ is the differential $1$-form on $N$ defined
by $\alpha_{\mid E_{\rm su}\oplus E_{\rm ss}}=0$ and $\alpha(Z)=1$,
called {\em Liouville's $1$-form}, then $\alpha\wedge
(d\alpha)^{n-1}$, where $n$ is the dimension of $M$, is a volume form
on $N$, which is invariant under the geodesic flow. Thus, the strong
stable leaves are contracted by the geodesic flow, and the strong
unstable leaves are dilated. See for instance \cite{KatHas95} for more
information.

\subsection{Discrete isometry groups}
\label{sec:discrisomgroups}

Let $\Ga$ be a discrete group of isometries of $\wt M$, which is {\em
  nonelementary}, that is, it preserves no set of one or two points in
$\wt M\cup\partial_\infty\wt M$.  To shorten the exposition, we will
assume in this survey that $\Ga$ has no torsion, though this
assumption is not necessary (see \cite{ParPau12JMD, ParPau13a,
  ParPau13b} for the extension), and is useful for some arithmetic
applications.

Let us denote the quotient space of $\wt M$ under $\Ga$ by $M=\Ga\bs
\wt M$, which is a smooth Riemannian manifold since $\Ga$ is torsion
free. We also say that the manifold $M$ is {\em nonelementary} if 
$\Ga$ is nonelementary.

\medskip We denote by $\Lambda\Ga$ the {\em limit set} of $\Ga$, that
is, the set of accumulation points in $\partial_\infty\wt M$ of any
orbit $\Ga x$ of a point $x$ of $\wt M$ under $\Ga$. It is the smallest
nonempty closed $\Ga$-invariant subset of $\partial_\infty\wt M$. The {\em
  critical exponent} of $\Ga$ is
$$
\delta_{\Ga}=\lim_{n\to+\infty}\frac 1n\ln\card\{\ga\in\Ga:
d(x_{0},\ga x_{0})\leq n\}. 
$$
The above limit exists (see \cite{Roblin02}), and the critical
exponent is positive (since $\Ga$ is nonelementary), finite (since
$M$ has a finite lower bound on its sectional curvatures, see for
instance \cite{Bowditch95}), and independent of the base point
$x_{0}$.

Since $\Ga$ acts without fixed points on $\wt M$, we have an
identification $\Ga\bs T^1\wt M=T^1M$, and we again denote by $\W^{\rm
  ss}, \W^{\rm s}, \W^{\rm su}$ and $\W^{\rm u}$ the continuous
foliations of $T^1M$ induced by the corresponding ones in $T^1\wt M$
(which are H\"older foliations when the sectional curvature of $M$ has
bounded derivatives).  We use the notation $(\flow{t})_{t\in\RR}$ also
for the geodesic flow on $T^1 M$.  We again denote by $\iota:T^1M\ra
T^1M$ the antipodal (flip) map $v\mapsto -v$, which also anti-commutes
with the geodesic flow.

Let us conclude this section by explaining some rigidity results on
the regularity of the foliations $\W^{\rm ss}, \W^{\rm s}, \W^{\rm
  su}$ and $\W^{\rm u}$.  Anosov has proved that if $M$ is compact,
then the vector subbundles $E_{\rm su}$ and $E_{\rm ss}$ are
H\"older-continuous (see for instance \cite[Th.~7.3]{PauPolSha} when
$M$ is only assumed to have pinched negative sectional curvature with
bounded derivatives). If $M$ is a compact surface, Hurder and Katok
\cite[Theo.~3.1, Coro.~3.5]{HurKat90} have proved that these
subbundles are $\operatorname{C}^{1,\alpha}$ for every $\alpha\in\;
]0,1[$ (see also \cite{HirPug75}), and that if they are
$\operatorname{C}^{1,1}$, then they are $\cinf$. Ghys
\cite[p.~267]{Ghys87} has proved that if $M$ is a compact surface, and
if the stable foliation of $T^1M$ is $\operatorname{C}^{2}$, then the
geodesic flow is $\cinf$-conjugated to the geodesic flow of a
hyperbolic surface. In higher dimension, we have the following result.

\btheo[Benoist-Foulon-Labourie \cite{BenFouLab92}] Let $M$ be a
compact negatively curved Riemannian manifold. If the stable foliation
of $T^1M$ is smooth, then the geodesic flow of $M$ is
$\cinf$-conjugated to the geodesic flow of a Riemannian symmetric
manifold with negative curvature. 
\cqfd 
\etheo

\section{Common perpendiculars of convex sets}
\label{sec:counting}

Let $M$ be a complete nonelementary connected Riemannian manifold of
dimension at least $2$, with pinched negative sectional curvature
$-b^2\le K\le -1$.  Let $\wt M\ra M$ be a universal Riemannian cover
of $M$, so that $\wt M$ is complete simply connected with the same
curvature bounds, and let $\Ga$ be its covering group, so that $\Ga$ is
a discrete, torsionfree, nonelementary group of isometries of $\wt M$.

\subsection{Convex subsets}
\label{subsec:convexsubset}

Let $\wt C$ be a nonempty {\it proper} (that is, different from $\wt
M$) closed convex subset of $\wt M$. Recall that a subset $A$ of $\wt
M$ is said to be {\em convex} if (the image of) any geodesic segment
with endpoints in $A$ is contained in $A$.  We denote by $\partial
\wt C$ the boundary of $\wt C$ in $\wt M$ and by $\partial_{\infty}\wt
C$ its set of points at infinity (the set of endpoints of geodesic
rays contained in $\wt C$). We say that the $\Ga$-orbit of $\wt C$ is
{\em locally finite} if, with $\Ga_{\wt C}$ the stabiliser of $\wt C$
in $\Ga$, for every compact subset $K$ of $\wt M$, the number of right
cosets $[\ga]\in\Ga/\Ga_{\wt C}$ such that $\ga\wt C$ meets $K$ is
finite.

Natural examples of convex subsets of $\wt M$ include the points, the
balls, the horoballs, the totally geodesic subspaces of $\wt M$ and
the convex hulls in $\wt M$ of the limit sets of nonelementary
subgroups of $\Ga$. Recall that the {\em convex hull} of a subset $A$
of $\partial_\infty\wt M$ with at least two points is the smallest
closed convex subset of $\wt M$ that contains $A$ in its set
of points at infinity.

Let $P_{\wt C}:\wt M\cup (\partial_{\infty} \wt M -\partial_{\infty}
\wt C) \to \wt C$ be the closest point map: if
$\xi\in \partial_{\infty} \wt M-\partial_{\infty}\wt C$, then $P_{\wt
  C} (\xi)$ is defined to be the unique point in $\wt C$ that
minimises the map $x\mapsto\beta_{\xi}(x,x_{0})$ from $\wt C$ to
$\RR$.  For every isometry $\ga$ of $\wt M$, we have $P_{\ga \wt
  C}\circ\ga=\ga\circ P_{\wt C}$.  The closest point map is continuous
in the topology of $\wt M\cup\partial_{\infty}\wt M$.

The {\it outer unit normal bundle} $\normal {\wt C}$ of $\wt C$ is
the subspace of $T^1\wt M$ consisting of the geodesic lines
$v\colon\RR\to \wt M$ with $v(0)\in\partial \wt C$, $v_+
\notin \partial_\infty \wt C$ and $P_{\wt C}(v_+) =v(0)$.  Note that
$\pi(\normal {\wt C})=\partial \wt C$, and that for all isometries
$\ga$ of $\wt M$, we have $\normal(\ga \wt C)=\ga\,\normal {\wt
  C}$. In particular, $\normal {\wt C}$ is invariant under the
isometries of $\wt M$ that preserve $\wt C$. When $\wt C=H\!B_-(v)$ is
the unstable horoball of $v\in T^1\wt M$, then $\normal{\wt C}$ is the
strong unstable manifold $W^{\rm su}(v)$ of $v$, and similarly,
$W^{\rm ss} (v)=\iota \,\normal H\!B_+(v)$.

The restriction of $P_{\wt C}$ to $\partial_{\infty} \wt M
-\partial_{\infty}\wt C$ (which is not necessarily injective) has a
natural lift to a homeomorphism
$$
P^+_{\wt C}:\partial_{\infty}\wt M-\partial_{\infty}\wt C\to
\normal {\wt C}
$$
such that $\pi\circ P^+_{\wt C}=P_{\wt C}$. The inverse of $P^+_{\wt
  C}$ is the {\it positive endpoint map} $v\mapsto v_+$ from $\normal
{\wt C}$ to $\partial_{\infty} \wt M - \partial_{\infty}\wt C$. In
particular, $\normal {\wt C}$ is a topological submanifold of $T^1\wt
M$, and a H\"older submanifold when the sectional curvature of $M$ has
bounded derivatives. For every $s\geq 0$, the geodesic flow induces a
homeomorphism $\flow s\colon \normal {\wt C}\to\normal\N_{s}\wt C$,
which is a H\"older homeomorphism when the sectional curvature of $M$
has bounded derivatives. For every isometry $\ga$ of $\wt M$, we have
$P^+_{\ga \wt C} \circ\ga= \ga\circ P^+_{\wt C}$.  When $\wt C$ has
nonempty interior and $\rm C^{1,1}$ boundary, then $\normal {\wt C}$
is the Lipschitz submanifold of $T^1\wt M$ consisting of the outer
unit normal vectors to $\partial \wt C$, and the map $P_{\wt C}$
itself is a homeomorphism (between $\partial_{\infty} \wt
M-\partial_{\infty}\wt C$ and $\partial \wt C$).  This holds when $\wt
C$ is the closed $\eta$-neighbourhood of any nonempty convex subset of
$\wt M$ with $\eta>0$ (see \cite[Theo.~4.8(9)]{Federer59},
\cite[p.~272]{Walter76}).

\medskip In this survey, we define a {\em properly immersed closed
  convex subset $C$ of $M$} as the data of a nonempty proper closed
convex subset $\wt C$ of $\wt M$, with locally finite $\Ga$-orbit, and
of the locally isometric proper immersion $C=\Ga_{\wt C}\bs\wt C\ra M$
induced by the inclusion of ${\wt C}$ in $\wt M$ and the Riemannian
covering map $\wt M\ra M$. (To simplify the exposition, we do not allow
in this survey the replacement of $\Ga_{\wt C}$ by one of its finite
index subgroups as it is done in \cite[\S 3.3]{ParPau13b}, even though this is
sometimes useful.) By abuse, when no confusion is possible, we will
again denote by $C$ the image of this immersion. We define $\normal
C=\Ga_{\wt C}\bs\normal {\wt C}$, which comes with a proper immersion
$\normal C\ra T^1M$ induced by the inclusion of $\normal {\wt C}$ in
$T^1\wt M$ and the covering map $T^1\wt M\ra T^1M$. By abuse also, we
will again denote by $\normal C$ the image of this immersion.

\subsection{The general counting problem}
\label{subsec:countingproblem}

Let $C_+,C_-$ be two properly immersed closed convex subsets of $M$.
A locally geodesic path $c:[0,T]\to M$ is a {\em common perpendicular}
from $C_-$ to $C_+$ if $\dot c(0)\in\normal C_-$ and $\dot c
(T)\in\iota\,\normal C_+$. For every $s\geq 0$, we denote by
$\Perp_{C_-,\,C_+}(s)$ the set of common perpendiculars from $C_-$ to
$C_+$ of length at most $s$. Each common perpendicular $c$ from $C_-$
to $C_+$ has a {\it multiplicity} $m(c)$, defined as follows. If $C_-$
and $C_+$ are the images in $M$ of two nonempty proper closed convex
subsets $\wt C_-$ and $\wt C_+$ of $\wt M$ with locally finite
$\Ga$-orbits, respectively, then $m(c)$ is the number of (left) orbits
under $\Ga$ of pairs $([\alpha],[\beta])$ in $\Ga/\Ga_{\wt C_-} \times
\,\Ga/\Ga_{\wt C_+}$ such that the closed convex subsets $\alpha \wt
C_-$ and $\beta\wt C_+$ have a (unique) common perpendicular whose image by $\wt M\ra M$ is $c$.  (Multiplicities are also
useful when $\Ga$ is allowed to have torsion, or when the stabilizers
$\Ga_{\wt C_\pm}$ are replaced by finite index subgroups, or when
$C_\pm$ is replaced by finite families of such convex subsets, see
\cite[\S 3.3]{ParPau13b} for a general version.)

In particular, any locally geodesic path is a common perpendicular of
its endpoints (with multiplicity $1$), since the outer unit normal
bundle of a point is equal to its unit tangent sphere. If $C_-$ and $C_+$
have nonempty interior and ${\rm C}^{1,1}$ smooth boundary (in the
appropriate sense for immersed subsets), the above definition of
common perpendicular agrees with the usual definition: a common
perpendicular exits $C_-$ perpendicularly to the boundary of $C_-$ at
its initial point and it enters $C_+$ perpendicularly to the boundary
of $C_+$ at its terminal point.

\medskip We study in this survey the asymptotic behaviour, as $s\to
+\infty$, of the number
$$
\N(s)=\N_{C_-,\,C_+}(s)=\sum_{c\in\Perp_{C_-,\,C_+}(s)} m(c)
$$ 
of common perpendiculars, counted with multiplicities, from $C_-$
to $C_+$, of length at most $s$. We refer to \cite[\S 3.3]{ParPau13b} for
more general counting functions.
 
Problems of this kind have been studied in various forms in the
literature since the 1950's and in a number of recent works, sometimes
in a different guise, as demonstrated in the examples below.  These
examples indicate that the general form of the counting results is
$\N(s)\sim \kappa\; e^{\delta s}$, where $\delta=\delta_{\Ga}$ is the
critical exponent of $\Ga$ and $\kappa>0$ is a constant.  Landau's
notation $f(s)\sim g(s)$ (as $s\to\infty$) means as usual that
$g(s)\ne 0$ for $s$ big enough, and that the ratio $\frac{f(s)}{g(s)}$
converges to $1$ as $s\to\infty$.

Observing that for $t\ge2\epsilon$, we have  
$$
\N_{\N_{\epsilon}(C_-),\,\N_{\epsilon}(C_+)}(t-2\epsilon)\leq
\N_{C_-,\,C_+}(t)\leq
\N_{\N_{\epsilon}(C_-),\,\N_{\epsilon}(C_+)}(t-2\epsilon)+
\N_{C_-,\,C_+}(2\epsilon)\;,
$$
we could replace the convex sets $C_-$ and $C_+$ by their
$\epsilon$-neighbourhoods for some fixed (small) positive $\epsilon$,
and then assume that $C_-$ and $C_+$ have $\operatorname{C}^{1,1}$
boundaries and use the more conventional definition of common
perpendicular. However, it is more natural to work directly with the
given convex sets instead of, for example, replacing points by small
balls.

\medskip Let $\Perp'(C_-,C_+)$ be the set $\Perp(C_-,C_+)$ where each
element $c$ has been replaced by $m(c)$ copies of it. The family
$(\ell(\alpha))_{\alpha\in\Perp'(C_-,C_+)}$ is called the {\it marked
  ortholength spectrum} of $(C_-,C_+)$. The set of lengths (with
multiplicities) of elements of $\Perp(C_-,C_+)$ is called the {\it
  ortholength spectrum} of $(C_-,C_+)$. This second set has been
introduced by Basmajian \cite{Basmajian93} (under the name ``full
orthogonal spectrum'') when $M$ has constant curvature, and $C_-$ and
$C_+$ are disjoint or equal, embedded, totally geodesic hypersurfaces or
embedded horospherical cusp neighbourhoods or embedded balls. Using
the complex lengths of the common perpendiculars between all closed
geodesics available in hyperbolic $3$-manifolds, and additional
combinatorial data, Meyerhoff \cite{Meyerhoff96} caracterizes the
isometry classes of closed hyperbolic $3$-manifolds. See also
\cite{Bridgeman11,BriKah10,Calegari11} for nice identities relating
the volume of $M$ and the ortholength spectrum, when $M$ is a compact
hyperbolic manifold with totally geodesic boundary and
$C_-=C_+=\partial M$. The aim of this paper is hence to survey the
asymptotic properties of the (marked) ortholength spectra.

\subsection{Counting orbit points in a ball}
\label{subsec:countorbitinball}

If $C_-=\{\bar p\}$ and $C_+=\{\bar q\}$ are singletons in $M$, then
$$
\N(s)=\card\big(B(p,s)\cap \Ga q\big)\;,
$$
for any lifts $p$ and $q$ of $\bar p$ and $\bar q$ in $\wt M$.  When
$\wt M=\hdr$ and $M$ is compact and orientable, Huber \cite[Satz
3]{Huber59} proved that
$$
\N(s)\sim\frac 1{4(g-1)}e^s,
$$ 
where $g$ is the genus of $M$.  His proof uses the Dirichlet series $
\sum_{\ga\in\Ga}\cosh^{-s}d(p,\ga q) $ and the Tauberian theorem of
Wiener-Ikehara \cite[\S 4,5]{Korevaar10}.

Margulis \cite[Theo.~2]{Margulis69} (see also \cite{Pollicott95})
generalised Huber's result for all compact connected negatively curved
manifolds of arbitrary dimension $n\geq 2$.  Note that $\delta=
\delta_\Ga$ is also the topological entropy of the geodesic flow of
$M$, since $M$ is compact. He showed that
$$
\N(s)\sim c(p,q)e^{\delta s}
$$ 
for some constant $c(p,q)$ which depends continuously on $p$ and
$q$. He proved that if $\wt M=\hnr$, then
\begin{equation}\label{eq:margulis}
\N(s)\sim\frac{\Vol(\SSS^{n-1})}{2^{n-1}(n-1)\Vol(M)}
\;e^{(n-1)s}\;.
\end{equation}
This agrees with Huber's result in dimension $n=2$ because the area of a
compact genus $g$ surface is $4\pi(g-1)$.  Margulis's proof of the
above result established the approach
\begin{center}
mixing $\to$ equidistribution $\to$ counting
\end{center} 
that has been used in most of the subsequent results.  Roblin
\cite[p.~56]{Roblin03} generalised Margulis's result when $\Ga$ is
nonelementary, the sectional curvature of $M$ is at most $-1$, and the
Bowen-Margulis measure of $T^1M$ is finite, and he has an expression
for the constant $c(p,q)$ in terms of the Patterson density and the
Bowen-Margulis measure, see \S \ref{sec:mainstatement} for more
details.  Lax and Phillips \cite{LaxPhi82} obtained an expression with
error bounds for the asymptotic behaviour of $\N(s)$ in terms of the
eigenvalues of the Laplacian on $\Ga\bs\hnr$ when $\Ga$ is {\it
  geometrically finite} (that is, in constant curvature, when $\Ga$ is
nonelementary and has a fundamental polyhedron with finitely many
sides).

\subsection{Counting common perpendiculars from a point 
to a totally geodesic submanifold}
\label{sec:pointandgeodesic}

Herrmann \cite[Theo.~I]{Herrmann62} proved for $\wt M=\hnr$, $M$
compact, $C_-=\{p\}$ a singleton, $C_+$ a compact totally geodesic
submanifold of dimension $k$, an asymptotic estimate
\begin{equation}\label{eq:herrmann}
\N(s)\sim\frac{2}{n-1}
\frac{\pi^{(n-k)/2}}{\Ga(\frac{n-k}2)}
\frac{\Vol(C_+)}{\Vol(M)}\,\frac{e^{(n-1)s}}{2^{n-1}}=
\frac{\Vol(\SSS^{n-k-1})\Vol(C_+)}{2^{n-1}\Vol(M)}\,\frac{e^{(n-1)s}}{n-1}\;,
\end{equation}
as $s\to+\infty$. Furthermore, he showed that the endpoints of the
common perpendiculars on the totally geodesic submanifold $C_+$ are
evenly distributed in terms of the Riemannian measure of $C_+$.  More
precisely, if $\Omega_+$ is a measurable subset of $C_+$ with boundary of (Lebesgue)
measure $0$, and if $\N_{p,\,\Omega_+}(s)$ is the number of those
common perpendiculars of $\{p\}$ and $C_+$ whose terminal endpoints are
contained in $\Omega_+$, then, as $s\to+\infty$,
\begin{equation}\label{eq:herrmann2}
\N_{p,\,\Omega_+}(s)\sim
\frac{\Vol(\SSS^{n-k-1})\Vol(\Omega_+)}
{2^{n-1}\Vol(M)}\,\frac{e^{(n-1)s}}{n-1}\,.
\end{equation}
The method of proof was a generalisation of that used by Huber.  The
asymptotic \eqref{eq:herrmann} was also treated in \cite{EskMcMul93}
as an illustration of their equidistribution result, see Theorem
\ref{theo:EM}.

Oh and Shah \cite{OhShaCircles} generalised Herrmann's result in
dimension $3$ for $\wt M=\htr$ and $\Ga$ geometrically finite (in which
case $\delta=\delta_\Ga$ is the Hausdorff dimension of the limit set
of $\Ga$, see for instance \cite{Bourdon95}). They showed that, as
$s\to+\infty$,
$$
\N(s)\sim c(\{p\},C_+)e^{\delta s}
$$ 
with a constant $c(\{p\},C_+)$ generalising that of Roblin's result
described above. Again, we postpone the description of the constant $
c(\{p\},C_+)$ until \S \ref{sec:mainstatement}.  This result is used
in \cite{OhShaCircles} to study $\Ga$-invariant families $(P_i)_{i\in
  I}$ of possibly intersecting circles in $\SSS^2$, called ``circle
packings'', that consist of a finite number of $\Ga$-orbits such that
the family of totally geodesic planes $(P^*_i)_{i\in I}$ in the ball
model of $\htr$ with $\partial_\infty P^*_i=P_i$ is locally finite.
They consider the counting function
$$
N(T)=\card\{i\in I\;:\;\curv(P^*_i)<T\},
$$
where the {\em spherical curvature} $\curv P^*_i$ is the cotangent of
the angle, at the origin $0$ of the ball model of $\htr$, between the
common perpendicular between $\{0\}$ and $P^*_i$, and any geodesic ray
from $0$ which is tangent to $P^*_i$ at infinity. Elementary
hyperbolic geometry (the angle of parallelism formula, see for
instance \cite[p.~147]{Beardon83}) implies that
$$
\curv P^*_i=\sinh d(0,P^*_i)\;,
$$
and thus, the above asymptotic estimate of $\N(s)$ is equivalent to
$N(T)\sim c(p,C_+) (2 T)^\delta$ as $T\to+\infty$.

\subsection{Counting common perpendiculars between 
horoballs}
\label{sec:horoballs}

When $C_-=C_+=\H$ is a {\em Margulis cusp neighbourhood} (that is, the
image by $\wt M\ra M$ of a $\Ga$-orbit of horoballs, centered at fixed
points of parabolic elements of $\Ga$, with pairwise disjoint
interiors), then results of \cite{BelHerPau01, HerPau04, Cosentino99,
  Roblin03} show that if $\Ga$ is geometrically finite (see
\cite{Bowditch95} for a definition in variable curvature), then
$$
\N(s)\sim c(\H)\,e^{\delta s}\;,
$$
as $t\to+\infty$, for some $c(\H)>0$.  Cosentino obtained explicit
expressions for the constant $c(\H)$ in special arithmetic cases:
$\Ga=\PSLZ$ acts on the upper halfplane model of $\wt M=\hdr$ and the
quotient space $\PSLZ\bs\hdr$ has a unique cusp that corresponds to
the orbit of $\infty$. The orbit of the subset $\wt \H$ of $\hdr$
consisting of points with vertical coordinates at least $1$ maps under
the quotient map to the maximal Margulis cusp neighbourhood of
$\PSLZ\bs\hdr$. The $\Ga$-orbit of $\wt \H$ consists of $\wt \H$ and
of the horoballs centered at all rational points $\frac pq$ with
$\gcd(p,q)=1$ and with Euclidean diameter $q^{-2}$.  The number of
such horoballs of diameter $q^{-2}$ modulo the stabiliser of $\infty$
(consisting of translations by the integers) is $\phi(q)$, where
$\phi$ is Euler's totient function. A classical result of Mertens on
the average order of $\phi$ (see for example
\cite[Theo.~330]{HarWri08}) implies that
$$
\N(s)=\frac 3{\pi^2}\,e^s+O(se^{s/2}),
$$
as $s\to+\infty$.  Similarly, if $\OOO_K$ is the ring of integers in
$K=\QQ(\sqrt d)$ with $d$ a negative squarefree integer and if $D_K$
is the discriminant of $K$, then $\Ga=\PSL(\OOO_K)$ acts on $\htr$ by
homographies as a cofinite volume discrete group of isometries. Let
$\H$ also be the image in $\Ga\bs \htr$ of the set of points in $\htr$
with vertical coordinates at least $1$. A generalisation of the above
argument gives, when $D_K\ne-3,-4$,
\begin{equation}\label{eq:cosentino}
\N(s)=\frac{\pi}{\sqrt{|D_K|}\,\zeta_K(2)}\,e^{2s}+O(e^{3t/2}),
\end{equation}
as $s\to+\infty$, where $\zeta_K$ is Dedekind's zeta function of
$K$. See \cite[Satz 2]{Grotz79} and \S 6.1, \S 6.2 in
\cite{Cosentino99} for a proof, and \cite{ParPauTou} for further
generalisations.

\subsection{Counting common perpendiculars of horoballs 
and totally geodesic submanifolds}
\label{sec:horoballgeodesic}

When $\wt M=\hnr$, $M$ has finite volume, $C_-$ is a Margulis cusp
neighbourhood of $M$ and $C_+$ is a finite volume totally geodesic
immersed submanifold of $M$ of dimension $k$ with $1\le k<n$, we
proved in \cite{ParPau12JMD} (see Theorem 1.1 and Lemma 3.1) the
following result, announced in \cite{ParPau10OR}.

\btheo[Parkkonen-Paulin \cite{ParPau12JMD}] As $s\to+\infty$,
\begin{align}
\N(s)&\sim\frac{\Vol(\SSS^{n-k-1})\Vol(C_-)\Vol(C_+)}
{\Vol(\SSS^{n-1})\Vol(M)} 
\;e^{(n-1)s}\nonumber
\\ &=
\frac{\Vol(\SSS^{n-k-1})\Vol(\partial C_-)\Vol(C_+)}
{\Vol(\SSS^{n-1})\Vol(M)} 
\;\frac{e^{(n-1)s}}{n-1}\;.\label{eq:parpau}\;\;\;\Box
\end{align}
\etheo

Oh and Shah \cite{OhSha12} studied a counting problem analogous to
the one described in \S \ref{sec:pointandgeodesic} for a bounded (to
simplify in this survey) family $\P$ of circles in $\RR^2$ that
consists of one (to simplify in this survey) orbit under a
nonelementary subgroup $\Ga$ of $\PSL_2(\CC)$ such that the family $\P^*$
of totally geodesic hyperplanes in the upper halfspace model of $\htr$
whose boundaries are the circles of $\P$ is locally finite. For any
circle $P\in\P$, let $\curv(P)$ be the reciprocal of the radius of
$P$, that is, the curvature of the circle $P$. For any $T>0$, let
$$
N(T)=\card\{P\in\P:\;\curv(P)<T\}\,.
$$
Oh and Shah showed \cite[Theo.~1.2]{OhShaCircles} that if $\delta>1$
and $\Ga$ is geometrically finite (see loc. cit.~for more general
assumptions)
\begin{equation}\label{eq:curvcount}
N(T)\sim c(\P)\, T^\delta
\end{equation}
for a constant $c(\P)>0$, as $T\to\infty$.  Furthermore, they proved
that the endpoints of the common perpendiculars are evenly distributed
on $\partial \wt C_-$ in the same sense as in Equation
\eqref{eq:herrmann2} in terms of a natural measure, which is the
skinning measure pushed to the boundary, see \S
\ref{sec:mainstatement}.

Let $\wt C_-$ be the horoball that consists of the points in the upper
half space model of $\htr$ whose vertical coordinates are at least $1$. Now,
$$
d(\wt C_-,\wt C_+)=\ln\curv(\wt C_+)
$$ 
for any hyperbolic plane $\wt C_+$ in $\P^*$. Hence, the result
\eqref{eq:curvcount} on circle packings has an interpretation as a
counting result for the common perpendiculars between the image $C_-$
of $\wt C_-$ and the image $C_+$ of the hyperplanes of $\P^*$ in
$M$. We then have
$$
N(T)=\N_{C_-,\,C_+}(\ln T)\,.
$$

\subsection{The density of integer points on homogeneous
  varieties}
\label{sec:DRS}

Let us denote a generic element of the Euclidean space $\RR^{n+1}$ by
$x=(x_{0},\bar x)$, where $x_{0}\in\RR$ and $\bar x=(x_1,\dots,
x_n)\in\RR^n$, and consider the quadratic form
$$
q(x)=-2x_{0}^2+\|x\|^2=-x_{0}^2+\|\bar x\|^2=-x_0^2+x_1^2+\dots+x_{n}^2
$$ 
of signature $(1,n)$.  The identity component $G=\operatorname{SO}_0
(1,n)$ of the special orthogonal group of the form $q$ is a connected
semisimple real Lie group with trivial center when $n\ge 2$.

Let $\RR^{1,n}=(\RR^{n+1},\langle\cdot,\cdot\rangle)$ be the
$(n+1)$-dimensional Minkowski space with the (indefinite) inner
product $$\langle x,y\rangle=-x_0y_0+\sum_{i=1}^nx_iy_i$$ where
$x=(x_0,x_1,\dots,x_n),y=(y_0,y_1,\dots,y_n)\in\RR^{n+1}$.  The
hyperboloid model of the $n$-dimensional real hyperbolic space $\hnr$
is the upper half $\{x\in\RR^{1,n}\;:\;q(x)=-1 ,\; x_0>0\}$ of the
hyperboloid with equation $q=-1$, endowed with the Riemannian metric
of constant sectional curvature $-1$ induced by the (positive
definite) restriction of the indefinite inner product $
\langle\cdot,\cdot\rangle$ to the tangent space of the hyperboloid.
The hyperbolic distance $d(x,y)$ of two points $x,y\in\hnr$ has a simple
expression in terms of the indefinite inner product: $\cosh
d(x,y)=-\langle x,y\rangle$.  The restriction to $\hnr$ of the (left)
linear action of $G$ on $\RR^{1,n}$ is the group of
orientation-preserving isometries of $\hnr$.

Oh and Shah proved the following counting result (and even a more
general one, see \cite[Th.~1.2]{OhSha13}) for linear orbits of
(nonelementary torsion free discrete) geometrically finite subgroups
$\Ga$ of $G$ on the level sets of the form $q$, generalising a special
case of a result of Duke, Rudnick and Sarnak \cite{DukRudSar93}: For
any $m\in\RR$, let $V_{m}= \{x\in\RR^{n+1}\;:\; q(x)=m\}$. Let $w\in
V_{m}-\{0\}$ be a vector such that the linear orbit $\Ga w$ is
discrete. If $\delta>1$, then by \cite[Coro.~1.6]{OhSha13}
\begin{equation}\label{eq:OSdensity}
\card\{y\in \Ga w:\|y\|<T\}\sim c(m) T^\delta\,,
\end{equation}
with a constant $c(m)>0$ similar to those in the previous cases.

The above counting result is equivalent to three results on counting
common perpendiculars, depending on the sign of $m$, as we will now
explain.  For convenience, we will restrict to the three essential
cases $m\in\{-1,0,1\}$.  Consider first the case $q(w)=-1$. Now, the
orbit of $w$ is contained in $\hnr$. For any $y\in\hnr$, we have
$\langle y,y\rangle=-y_{0}^2+\|\bar y\|^2=-1$. Thus, $\|y\|^2=
2y_0^2+1$, and we have $\cosh d(y,(1,0))=-\langle (1,0),y\rangle=
y_{0} \sim\frac 1{\sqrt 2}\|y\|$ as $\|y\|\to+\infty$. Therefore, the
asymptotic \eqref{eq:OSdensity} gives an asymptotic count of orbit
points as in the results of Margulis and Roblin, see \S
\ref{subsec:countorbitinball}.

If $q(w)=0$, then $w$ lies in the light cone of $q$ and it defines a
horosphere 
$$
H_{w}=\{y\in\hnr:\langle y,w\rangle=-1\}\;.
$$ 
Using a rotation with fixed point at $(0,1)$, we can assume that
$w=(w_{0},w_{0},0)\in\RR\times\RR\times\RR^{n-1}$, with $w_0\neq
0$. Now
$$
H_{w}=\{y\in\hnr:y_{1}=y_{0}-\frac
1{w_{0}}\}=\{y\in\RR^{1,n}:y_{0}=\frac 12(w_{0}+\frac 1{w_{0}}),\,
y_{1}=\frac 12(w_{0}-\frac 1{w_{0}}) \}\,. 
$$
By the symmetry of the situation, it is clear that 
\begin{align*}
d((1,0),H_{w})= &d((1,0),(\frac 12(w_{0}+\frac 1{w_{0}}),\frac
12(w_{0}-\frac 1{w_{0}}),0))=\operatorname{arcosh}(\frac
12(w_{0}+\frac 1{w_{0}}))\\ 
= & \ln w_{0}=\ln\|w\|-\ln 2\,,
\end{align*}
so the asymptotic for the norms of points in the orbit of a point in
the light cone is equivalent to an asymptotic of the distance of an
orbit of horoballs from a point. The same counting problem is also
considered by Kontorovich and Oh in \cite{KonOh11}, and earlier in
\cite{Kontorovich09} in the two-dimensional case.

In the third case, when $q(w)=1$, the vector $w$ defines a totally
geodesic hyperplane $w^\perp=\{y\in\hnr:\langle y,w\rangle=0\}$ in
$\hnr$. As in the two cases above, one can check that this asymptotic
is equivalent to the asymptotic count of geodesic arcs starting at a
fixed point and ending perpendicularly at an orbit of totally geodesic
hyperplanes.

\section{Using Eskin-McMullen's equidistribution theorem}
\label{sec:eskinmcmullen}

In order to prove the kind of asymptotic results described in \S
\ref{sec:counting}, following Margulis \cite{Margulis69,Margulis04},
one usually proves first an appropriate equidistribution result using
mixing, and this result is then used to study the common
perpendiculars.

Eskin and McMullen \cite[Th.~1.2]{EskMcMul93} proved a very general
equidistribution theorem for Lie groups orbits using mixing properties
and a technical ``wave front lemma'' in affine symmetric spaces.

\btheo[Eskin-McMullen]
\label{theo:EM} Let $G$ be a connected semisimple real Lie group with
finite center. Let $\sigma: G\to G$ be an involutive Lie group
automorphism, and $H$ be its fixed subgroup.  Let $\Ga$ be a lattice in
$G$ and let $m$ be the unique $G$-invariant probability measure on
$\Ga\bs G$. Assume that the projection of $\Ga$ to $G/G'$ is dense for
all noncompact connected normal Lie subgroups $G'$ of $G$, and that
$\Ga\cap H$ is a lattice in $H$. Let $Y=(\Ga\cap H)\bs H$ and let
$\mu_{g}$ be the image by the right multiplication by $g$ of the
unique $H$-invariant probability measure on $Y$. Then, for every $f:
\Ga\bs G\to\RR$ which is continuous with compact support,
$$
\int_{Yg}f(h)\;d\mu_{g}(h)\to\int_{\Ga\bs G}f(x)\;dm(x),
$$
as $g$ goes to infinity in $H\bs G$. \cqfd
\etheo

\noindent This result is used in \cite{EskMcMul93} to prove a result
of Duke, Rudnick and Sarnak \cite{DukRudSar93} on counting integral
points on homogeneous varieties, see \S \ref{sec:DRS}.

In \cite{ParPau12JMD}, we proved the following equidistribution result
using mixing and hyperbolic geometry, as a tool to prove the
asymptotic estimate \eqref{eq:parpau}. A modification of the proof in
\cite{ParPau12JMD} enabled us to prove the general equidistribution
result in variable curvature in \cite{ParPau13a}, whose tools are also
used for the counting result of \cite{ParPau13b} that we describe in
\S \ref{sec:mainstatement}. Here, at the instigation of Hee Oh,
we present a different proof using Theorem \ref{theo:EM}, which also
serves as an illustration of the use of Theorem \ref{theo:EM}.

\btheo[Parkkonen-Paulin \cite{ParPau12JMD}] 
\label{theo:equidistribution} Let $M$ be a complete connected
hyperbolic ma\-nifold with finite volume.  Let $C$ be a nonempty
proper totally geodesic immersed submanifold of $M$ with finite
volume. The induced Riemannian measure on $g^t\normal C$
equidistributes to the Liouville measure as $t\ra+\infty$:
$$
\Vol_{g^t\normal C}/\|\Vol_{g^t\normal C}\|\;\;
\stackrel{*}{\rightharpoonup}\;\;\Vol_{T^1M}/\|\Vol_{T^1M}\|\;. 
$$
\etheo

\noindent More general versions of the above result appear in
\cite[Theo.~1.8]{OhSha13}, \cite[Theo.~17]{ParPau13a} and 
\cite[Theo.~20]{ParPau13b} in the presence of potentials.

We use below the notation introduced in \S \ref{sec:DRS}. Before
giving the proof of this result, let us review some preliminaries on
the action on $T^1\hnr$ of the orientation-preserving isometry group
$G= \operatorname{SO}_0(1,n)$ of the hyperboloid model of $\hnr$,
where $n\geq 2$.  Let $(e_0,e_1,\dots,e_n)$ be the canonical basis of
$\RR^{1,n}$, and let $w_0=(1,0,\dots,0)\in\hnr$.  For any $1\le k< n$,
we embed $\hkr$ isometrically in $\hnr$ as the intersection of $\hnr$
with the linear subspace given by the equations
$x_{k+1}=x_{k+2}=\dots=x_n=0$.  For any $p\in\NN$, let $I_p$ be the
$p\times p$ identity matrix. Let $H_k$ be the subgroup of $G$ that
consists of the fixed points of the involution $\sigma_k\colon G\to G$
defined by $\sigma_k(g)= J_kgJ_k^{-1}$, where
$J_k=\begin{pmatrix}I_{k+1}& 0\\ 0& -I_{n-k}
\end{pmatrix}$. Note that $H_k$ is isomorphic to
$(\operatorname{O}(1,k) \times \operatorname{O}(n-k))\cap G$, hence
contains $\operatorname{SO}_0(1,k) \times \operatorname{SO}(n-k)$ with
index $2$. Let us identify $\operatorname{SO}(n-1)$ with its image in
$\operatorname{SO}(n)$, and similarly $\operatorname{SO}(n)$ with its
image in $G$, by the maps $x\mapsto\begin{pmatrix}1& 0\\ 0& x
\end{pmatrix}$. Let  $\lambda_G$ and $\lambda_{H_k}$ be fixed left
Haar measures on  $G$ and $H_k$.

The group $G$ acts transitively on $T^1\hnr$ and its action commutes
with the geodesic flow, the stabiliser of $e_1\in T^1\hnr$ being
$\operatorname{SO}(n-1)$.  Note that $H_k$ is the subgroup of $G$
which preserves $\hkr$. It acts transitively on the unit normal bundle
$\normal\hkr$.

The orbital map $g\mapsto ge_1$ from $G$ to $T^1\hnr$ induces a
diffeomorphism $\overline\varphi\colon G/\operatorname{SO}(n-1)\ra
T^1\hnr$ which is equivariant for the left actions of $G$.  The
commutativity of the diagram
$$
\begin{array}{ccc} 
G/\operatorname{SO}(n-1) &\longrightarrow & G/\operatorname{SO}(n)
\smallskip \\ 
\;\;\downarrow\simeq\overline\varphi  & & \;\downarrow\simeq 
\smallskip \\
T^1\hnr &\longrightarrow & \hnr
\end{array}
$$
and the fact that the Riemannian measure of $\SSS^{n-1}\simeq
\operatorname{SO}(n)/\operatorname{SO}(n-1)$ is the unique (up to
multiplication by a positive constant) positive Borel measure which is
invariant under rotations imply that the image of $\lambda_G$ by the
smooth map $g\mapsto \bar\varphi(g\operatorname{SO}(n-1))$ is a
multiple of $\Vol_{T^1\hnr}$.

Consider the one-parameter subgroup $(a_t)_{t\in\RR}$ of $G$, where
$$
a_t=\begin{pmatrix}\cosh t& \sinh t & 0\\ \sinh t&\cosh t& 0\\
  0 & 0 & I_{n-1}\end{pmatrix}\,.
$$  
The action $g\mapsto ga_t$ of $a_{t}$ by right translations on $G$
commutes with that of $\operatorname{SO}(n-1)$.  A calculation in
hyperbolic geometry shows that $a_te_1$ is the image of the unit
tangent vector $e_1\in T^1_{w_0}\hnr$ under the geodesic flow
$\flow{t}$ in $T^1\hnr$. Thus, by equivariance, $\overline \varphi
(ga_t \operatorname{SO}(n-1))= \flow{t}\overline \varphi(g
\operatorname{SO}(n-1))$ for all $g\in G$ and all $t\in\RR$.  Let us
fix a group element $r\in\operatorname{SO}(n)$ which maps $e_1$ to
$e_{n}$. As the measures under consideration are induced by
differential forms, homogeneity arguments imply that the image measure
of $\lambda_{H_{k}}$ by the smooth $H_{k}$-equivariant map $h\mapsto
\overline\varphi (hra_t\operatorname{SO}(n-1))$ is a multiple of
$\Vol_{g^t\normal\hkr}$.

\medskip
\noindent{\bf Proof of Theorem \ref{theo:equidistribution}.} By
additivity, we can assume that $C$ is connected.  Let $\overline{M}\ra
M$ be the Riemannian orientation cover of $M$ (which is the identity
map if $M$ is orientable), and let $\overline C\ra C$ be the one of
$C$, so that $\overline C $ is a connected immersed totally geodesic
submanifold of $\overline M$. As the image measures by the finite
cover $T^1\overline M\ra T^1M$ of $\Vol_{g^t\partial^1_+\overline{C}}$
and $\Vol_{T^1\overline{M}}$ are $\Vol_{g^t\partial^1_+C}$ and
$\Vol_{T^1M}$, respectively, we only have to show that the Riemannian
measure of $g^t\partial^1_+\overline{C}$ equidistributes to the
Liouville measure of $T^1\overline{M}$ as $t\ra +\infty$. We may,
therefore, assume that $M$ and $C$ are oriented.

Let us fix a universal Riemannian cover $\hnr\ra M$. Its covering
group $\Ga$ is a lattice in $G$, since $M$ has finite volume. We may
assume that the image of $\hkr$ under this covering map is $C$. We
define $H=H_{k}$ as above. Since $C$ has finite volume, $\Gamma\cap H$
is a lattice in $H$.  Since $r$ fixes $e_0$ and sends $e_1$ to $e_n$
which is perpendicular to $\hkr$, the map $t\mapsto \pi(ra_te_1)$ from
$[0,+\infty[$ to $\hnr$ is a geodesic ray starting perpendicularly to
$\hkr$. Since $H$ is the stabiliser of $\hkr$, the map $t\mapsto
Hra_te_1$ from $[0,+\infty[$ to $H\bs\hnr$ tends to infinity. Hence
the map $t\mapsto Hra_t$ from $[0,+\infty[$ to $H\backslash G$ tends
to infinity.

Since the connected semi-simple real Lie group $G$ has trivial center,
and only one noncompact factor, the projection of $\Ga$ to $G/G'$ is
dense for all noncompact connected normal Lie subgroups $G'$ of
$G$. (In fact, $G$ has no compact factor, and any lattice in $G$ is
irreducible, see for instance \cite{Mostow73}).

We can now use Theorem \ref{theo:EM} to conclude that as $t$ tends to
$+\infty$, the measure $\mu_t$ on $\Gamma\backslash G$ with support
$\Gamma H r a_t$ which is defined to be the translate by $r a_{t}$ of
the unique $H$-invariant probability measure on $(\Ga\cap H)\bs H$,
equidistributes towards the probability measure $m$ on $\Ga \bs G$
induced by $\lambda_G$. Let $p:\Gamma\backslash G\ra \Gamma\backslash
G/\operatorname{SO}(n-1)$ be the canonical projection.  The
$G$-equivariant diffeomorphism $\overline\varphi: G/\operatorname{SO}
(n-1)\ra T^1\HH^n_\RR$ induces a diffeomorphism $\varphi: \Ga\bs
G/\operatorname{SO}(n-1)\ra T^1M$ such that $\varphi (\Ga g a_t
\operatorname{SO}(n-1))=g^t\varphi(\Ga g \operatorname{SO}(n-1))$ for
all $g\in G$ and $t\in\RR$.  By the homogeneity argument just before
the beginning of the proof and covering arguments, and since direct
images of measures preserve the total masses, we have
$\varphi_*(p_*\mu_t)= \frac{1}{\Vol(g^t\partial^1_+C)}
\Vol_{g^t\partial^1_+C}$ and $\varphi_* (p_*m)= \frac{1}{\Vol(T^1M)}
\Vol_{T^1M}$. As taking direct images of measures by a given
continuous map is continuous in the weak-* topology, the measures
$\frac{1}{\Vol(g^t\partial^1_+C)} \Vol_{g^t\normal C} = (\varphi\circ
p)_*\mu_t$ equidistribute to $(\varphi\circ p)_*m =
\frac{1}{\Vol(T^1M)} \Vol_{T^1M}$ in $T^1M$ as $t$ tends to $+\infty$,
which is what we wanted to prove.  \cqfd

\section{Arithmetic applications} 
\label{sec:arithmapplic}

If the manifold $M$ is arithmetically defined, many counting results
for common perpendiculars have an arithmetic interpretation. In
this section, we will review some of these arithmetic
applications. The arithmetically defined groups in this section will,
in general, have torsion. This is not a problem however, as the
geometric counting result used in the various cases below is indeed
valid in this more general context.  In certain cases, the interaction
has also worked in the opposite direction, as evidenced by the results
of Cosentino in \S \ref{sec:horoballs}.

\subsection{Counting representations of  integers by 
binary quadratic forms}
\label{sec:quadratic}

Let $Q(X,Y)=aX^2+bXY+cY^2$ be an integral binary quadratic form with
{\em discriminant} $\Delta=b^2-4ac$. An element $x\in\ZZ^2$ is a {\em
  representation} of an integer $n$ by $Q$ if $Q(x)=n$, and the
representation is {\em primitive} if the components of $x$ are
relatively prime. If $Q$ is positive definite (equivalently, if
$\Delta<0$ and $a>0$), then the number $N(t)$ of representations of integers
that are at most $t$ equals the number of lattice points of $\ZZ^2$
inside the ellipse defined by the equation $Q(x)=t$.  The asymptotics
of this number (Gauss' circle problem) have been studied extensively,
and the best known result with an error bound
$$
N(t)=\frac{2\pi}{\sqrt{-\Delta}}\;t+O(t^{131/416})
$$ 
is due to Huxley \cite{Huxley03}. Gauss already had a solution with a
worse bound on the error term.

The modular group $\SLZ$ acts on the right by precomposition on the
set of binary quadratic forms, preserving the discriminant, and
linearly on the left on $\ZZ^2$.  Let us assume that $Q$ is {\em
  primitive} (that is, the coefficients $a$, $b$ and $c$ are
relatively prime), {\it indefinite} (that is, $\Delta>0$) and not the
product of two integral linear forms.  The stabiliser of a form $f$ in
$\SLZ$, called the {\em group of automorphs} of $f$, is
\begin{align*}
\operatorname{SO}(Q,\ZZ)= & \{\ga\in\SLZ:Q\circ\ga=Q\}\\
= & \Big\{\ga_{Q,t,u}= 
\begin{pmatrix}
\dfrac{t-bu}2 & -cu\\ au & \dfrac{t+bu}2
\end{pmatrix}:\;t,u\in\ZZ,\;\; t^2-\discr u^2=4\Big\}\;,
\end{align*}
see for instance \cite[Theo.~202]{Landau58}. This group is infinite
and thus any nonzero integer that is represented by the form $Q$ is
represented infinitely many times.  Accordingly, in the generalisation
of the circle problem for these forms, one counts the number of orbits
of lattice points under the linear action of $\operatorname{SO}(Q,
\ZZ)$ between the hyperbolas defined by the equations $|Q(x)|=t$. With
$\P$ the set of relatively prime elements of $\ZZ^2$, let
$$
\tilde\Psi_Q(t)=\card\big(\operatorname{SO}(Q,\ZZ)\bs
\{x\in\ZZ^2\;:\; |Q(x)|\le t\}\big)
$$
and
$$
\Psi_Q(t)=\card\big(\operatorname{SO}(Q,\ZZ)\bs
\{x\in\P\;:\; |Q(x)|\le t\}\big)
$$
be the counting functions of all the representations and of the
primitive representations by $Q$. The asymptotics of $\tilde\Psi_Q(t)$
are also known, see for example \cite[p.~164]{Cohn62} for a proof.
It turns out that the asymptotic result on the counting function
$\N(s)$ for a horoball and a totally geodesic subspace can be used to
give a different proof of this result.

Let us first observe that an asymptotic result for the primitive
representations implies one for all representations. Assume that
$\Psi_Q(t)=ct+O(t^{1-\epsilon})$ for some $c,\epsilon>0$.  For any
$k\in\NN$, let
$$
\Psi_{Q,k}(t)=\card\big(\operatorname{SO}(Q,\ZZ)\bs
\{x=(x_1,x_2)\in\ZZ\;:\;
\gcd(x_1,x_2)=k,\,|Q(x)|\le t\}\big)\,.
$$
Now, $\Psi_{Q,k}(t)=\Psi_{Q}(k^{-2}t)$ and
$$
\tilde\Psi_Q(t)=\sum_{k=1}^\infty \Psi_{Q,k}(t)=\sum_{k=1}^\infty
\Psi_{Q}(k^{-2}t)=\sum_{k=1}^\infty  ck^{-2}t+ O((k^{-2}t)^{1-\epsilon})
=c\,\zeta(2)t+O(t^{1-\epsilon}).
$$

We will now explain how to obtain an asymptotic estimate for
$\Psi_Q(t)$ from the solution of the geometric counting problem of \S
\ref{sec:horoballgeodesic}.  We use the upper halfplane model of
$\hdr$. The subgroup $\operatorname{PSO} (Q, \ZZ)$ of $\PSLZ$ is a
cyclic group generated by a hyperbolic element. Its index $i_{Q}$ in
the stabiliser $\Ga_{Q}$ of the geodesic line $C_{Q}$ invariant under
$\operatorname{PSO}(Q,\ZZ)$ is either $2$ or $1$ depending on whether
or not the corresponding locally geodesic line on the modular surface
$\PSLZ\bs \hdr$ passes through the cone point of order $2$.  Let
$(t_{Q},u_{Q})$ be the fundamental solution of the Pell-Fermat
equation $t^2-\discr u^2=4$, and let $R_Q=\ln \frac{t_{Q} +
  u_{Q}\sqrt\discr}2$ be the {\em regulator} of $Q$.  It is easy to
check that the length of the closed geodesic $\Ga_{Q}\bs C_{Q}$ is
$\frac{2R_{Q}}{i_{Q}}$.

The stabiliser $U$ of $(1,0)\in\ZZ^2$ for the linear action of $\SLZ$
is the subgroup that consists of integral upper triangular unipotent
matrices. Geometrically, the image $\Ga_{\infty}$ of $U$ in $\PSLZ$ is
the stabiliser in $\PSLZ$ of the horoball $\H=\{z\in\CC\;:\;\Im\;
z\ge1\}$ in $\hdr$.  The horoball $\H$ is {\em precisely invariant},
that is, each element of $\PSLZ$ either preserves $\H$ or maps $\H$ to
a horoball whose interior is disjoint from $\H$. As $\H$ is the
maximal such horoball at $\infty$, it corresponds to the maximal
Margulis cusp neighbourhood of the unique cusp of the quotient space
$\PSLZ\bs\hdr$.

The length of the common perpendicular of $\H$ and the hyperbolic line
$C_{Q}$ stabilised by $\operatorname{SO}(Q,\ZZ)$ is $\ln\frac{2|a|}
{\sqrt\Delta}$. For all $\ga=\pm\Big(\!\!\begin{array}{cc}A & B \\
  C & D\end{array}\!\!\Big)$ in $\operatorname{PSL}_2(\ZZ)$, a simple
computation (see \cite[Lem.~4.2]{ParPau12JMD}) shows that the length
of the common perpendicular of $\H$ and the image under $\ga$ of
$C_{Q}$ is
$$
\ln\frac{2}{\sqrt{\discr}}\;|Q(D,-C)|\;.
$$
Thus, Corollaire 3.9 of \cite{ParPau12JMD}, which generalises the
result of Equation \eqref{eq:parpau} to the case of groups with
torsion, and \cite[\S 6]{ParPau13b} (see \S \ref{sec:errorterm}) for
the error term, give that there exists $\kappa>0$ such that
\begin{align*}
\Psi_Q(s)&= i_Q\operatorname{Card}\big\{[\ga]\in\Ga_{\infty}\backslash
\Ga/\Ga_{Q}\;:\; \;d(\H_\infty, \ga C_Q)\leq \ln
\Big(\frac{2}{\sqrt{\discr}}\;s\Big)\big\}\\
&\sim i_Q\;\frac{\Vol(\SSS^{0})
\Vol(\Ga_{\infty}\backslash\H_\infty)\Vol(\Ga_{Q}\backslash C_Q)}
{\Vol(\SSS^{1})\Vol(\Ga\backslash\HH^2_\RR)} 
\;\Big(\frac{2}{\sqrt{\discr}}\;s\,(1+\operatorname{O}(s^{-\kappa}))\Big)
\\ &=\frac{12\, R_{Q}}{\pi^2\sqrt{\discr}}\; 
s\,(1+\operatorname{O}(s^{-\kappa}))\,,
\end{align*}
see \cite[\S 4]{ParPau12JMD} for more details and more applications,
in particular to counting representations satisfying congruence
relations, and to \cite[\S 6]{ParPau13b} for error terms.

\subsection{Counting representations of integers by binary  
Hermitian  forms} 
\label{sec:hermitian}

A function $f:\CC^2\ra\RR$ is a {\em binary Hermitian form} if there
are constants $a,c\in\RR$ and $b\in\CC$, called the {\em coefficients}
of $f$, such that for all $u,v\in\CC$,
\begin{equation}\label{eq:form}
f(u,v)=a|u|^2+2\,\Re(bu\ov v) +c|v|^2\,=
\begin{pmatrix}\bar u & \bar v \end{pmatrix}
\begin{pmatrix} a& b\\ \bar b & c \end{pmatrix}
\begin{pmatrix}u \\  v \end{pmatrix}.
\end{equation}
Let $K$ be an imaginary quadratic number field, with discriminant
$D_K$ and ring of integers $\OOO_K$. If the coefficients of the
Hermitian form $f$ satisfy $a,c\in\RR\cap\OOO_{K}=\ZZ$ and $b\in
\OOO_{K}$, then we say that $f$ is {\em integral} (over $K$). It is
easy to check that the values of the restriction of an integral binary
Hermitian form to $\OOO_{K}\times\OOO_{K}$ are rational integers.  If
the {\em discriminant} $\Delta(f)=|b|^2-ac$ of $f$ is positive, then
we say that $f$ is {\em indefinite}, which is equivalent to saying
that $f$ takes both positive and negative values.

A binary Hermitian form naturally gives rise to a quaternary quadratic
form. The representations of integers by positive definite quaternary
quadratic forms have been studied for a long time (including
Lagrange's four square theorem, see also the work of Ramanujan as in
\cite{Kloosterman27}).

In the case of indefinite forms, the counting problem is again
complicated by the presence of an infinite group of automorphs: The
group $\SLOK$ acts on the right by precomposition on the set of
(indefinite) integral binary Hermitian forms, and the stabiliser of
such a form under this action is, analogously to the case of binary
quadratic forms treated in \S \ref{sec:quadratic}, called the {\it
  group of automorphs} of the form and denoted by $\operatorname{SU}_f
(\OOO_K)$.  The Bianchi group $\PSLOK$ acts discretely on the upper
halfspace model of $\htr$, with finite covolume. Now, the image in
$\PSLOK$ of the group of automorphs of a fixed indefinite integral
binary Hermitian form $f$ is a Fuchsian subgroup that preserves a real
hyperbolic plane $\C(f)$ whose boundary at infinity is the circle
$$
\C_\infty(f)=\{[u:v]\in\PP^1(\CC)=\partial_{\infty}\htr\;:\;f(u,v)=0\}\,.
$$
The group of automorphs $\operatorname{SU}_f(\OOO_K)$ is an arithmetic
group that acts on $\C(f)$ with finite covolume.

The cusps of a Bianchi group $\PSLOK$ are in a natural bijective
correspondence with the ideal classes of $K$, see for example Theorem
2.4 in Chapter 7 of \cite{ElsGruMen98}.  Let $x,y\in\OOO_{K}$ be not
both zero, so that $[x:y]\in\PP^1(\CC)$ is a cusp of $\PSLOK$. Then,
if $y=0$, the horoball $\H$ that consists of those points in the upper
halfspace model $\htr$ whose vertical coordinate is at least $1$ is
precisely invariant, and if $y\ne0$, then there is some $\tau>0$ such
that the horoball $\H$ centered at $\frac{x}{y}$ of Euclidean height $
\tau$ is precisely invariant. Analogously with the case of indefinite
binary quadratic forms, for any $g\in\SLOK$, the hyperbolic distance
between $\H$ and $\C_\infty(f\circ g)=g^{-1}\C_\infty(f)$ is $\ln
\frac{|f\circ g(x,y)|}{\tau|y|^2\sqrt{\Delta(f)}}$, and, as in the
case of binary quadratic forms, we find a connection between
representing integers by $f$ and the counting problem of \S
\ref{sec:horoballgeodesic}.

We define a counting function of the representation of integers for
each nonzero fractional ideal $\mmm$ of $K$.  For every $u,v\in K$,
let $\langle u,v\rangle$ be the $\OOO_K$-module they generate.  For
every $s>0$, we consider the integer
$$
\psi_{f,\mmm}(s)=\card\;\;_{\operatorname{SU}_f(\OOO_K)}\bs
\big\{(u,v)\in\mmm\times\mmm\;:\;(N\mmm)^{-1}|f(u,v)|\leq s,
\;\;\;\langle u,v\rangle=\mmm\big\}\;.
$$
Generalising the argument used for binary quadratic forms (see \S
\ref{sec:quadratic}), we can again use the generalisation of Equation
\eqref{eq:parpau} to obtain an asymptotic expression for
$\psi_{f,\mmm}(s)$.

\btheo[Parkkonen-Paulin \cite{ParPau11BLMS}\cite{ParPau13b}] 
\label{theo:mainHerm}
There exists $\kappa>0$ such that, as $s$ tends to $+\infty$,
$$
\psi_{f,\mmm}(s)\;=\; \frac{\pi\;\covol(\operatorname{SU}_f(\OOO_K))}
{2\;|D_K|\;\zeta_K(2)\;\Delta(f)}\;\;\;s^2
\,(1+\operatorname{O}(s^{-\kappa}))\;.\;\;\;\Box
$$
\etheo

Here $\covol(\operatorname{SU}_f(\OOO_K))$ is the area of the quotient
of the hyperbolic plane $\C(f)$ in $\htr$ by the group of automorphs
of $f$, and $\zeta_K$ is Dedekind's zeta function of $K$.  In the proof,
after applying Equation \eqref{eq:parpau}, we use the fact that there
is an explicit formula (essentially due to Humbert) for the volume
$$
\Vol(\,\PSLOK\bs\htr)=\frac{1}{4\pi^2}|D_K|^{3/2}\zeta_K(2)\,.
$$ 
See \cite{Sarnak83} for a proof of this formula using Eisenstein
series, and \S 8.8 and \S 9.6 of \cite{ElsGruMen98} for further
proofs.  The following corollary follows immediately by taking
$\mmm=\OOO_K$: If $\P_K$ is the set of relatively prime pairs of
integers of $K$, then
$$
\card\;\;_{\operatorname{SU}_f(\OOO_K)}\bs
\big\{(u,v)\in\P_K\;:\;|f(u,v)|\leq s\big\}
\;=\; \frac{\pi\;\covol(\operatorname{SU}_f(\OOO_K))}
{2\;|D_K|\;\zeta_K(2)\;\Delta(f)}\;\;\;
s^2\,(1+\operatorname{O}(s^{-\kappa}))\;,
$$
as $s$ tends to $+\infty$.

In general, one could compute the covolume of the group of automorphs
$\operatorname{SU}_f(\OOO_K)$ with the aid of Prasad's formula in
\cite{Prasad89}.  Maclachlan and Reid \cite{MacRei91} computed the
covolumes of all stabilisers in $\operatorname{PSL}(\QQ(i))$ of
Euclidean halfspheres in the upper halfspace model of $\htr$ centered
at $0$ with Euclidean radius $\sqrt D$, where $D$ is a rational
integer.  This result can be used to obtain an even more explicit
expression of the asymptotic formula of Theorem \ref{theo:mainHerm}
when $K= \QQ(i)$: A constant $\const(f)\in\{1,2,3,6\}$ is defined as
follows.  If $\discr(f)\equiv 0\mod 4$, let $\const(f)=2$. If the
coefficients $a$ and $c$ of the form $f$ as in Equation
\eqref{eq:form} are both even, let $\const(f)=3$ if $\discr(f)\equiv
1\mod 4$, and let $\const(f)$ be the remainder modulo $8$ of
$\discr(f)$ if $\discr(f)\equiv 2 \mod 4$.  In all other cases, let
$\const(f)=1$. The class number of $\QQ(i)$ is $1$, and there is just
one counting function to be considered.  We prove in
\cite[Coro.~3]{ParPau11BLMS}, and \cite[\S 6]{ParPau13b} (see \S
\ref{sec:errorterm}) for the error term that, if $K=\QQ(i)$, there
exists $\kappa>0$ such that, as $s$ tends to $+\infty$,
$$
\psi_{f,\,\ZZ[i]}(s)=\frac{\pi^2}{8\;\const(f)\;\zeta_{\QQ(i)}(2)}\,
\prod_{p|\discr(f)} \big(1+\bigg(\frac{-1}p\bigg)p^{-1}\big)\;\;
s^2\,(1+\operatorname{O}(s^{-\kappa}))\;.
$$
Here $p$ ranges over the odd positive rational primes and
$\big(\frac{-1}p\big)$ is the Legendre symbol of $-1$ modulo $p$.  We
refer to \cite{ParPau11BLMS} for more details and more applications,
including counting representations satisfying congruence conditions,
and to \cite[\S 6]{ParPau13b} for error terms.

\subsection{Counting quadratic irrational in orbits of modular 
groups}
\label{subsec:fibonacci}

The group $\PSLZ$ acts transitively on the rational real numbers, but
not transitively on the irrational algebraic real numbers of a given
degree. Hence, counting results (for appropriate complexities) of
algebraic irrationals within an orbit of $\PSLZ$ is an interesting
problem, and we give some solutions in \cite{ParPau12JMD} in the
quadratic case. Similar problems occur for quadratic irrational
complex numbers under the action of (congruence subgroups of) Bianchi
groups, and we illustrate them by the following result.

Let $\phi=\frac{1+\sqrt{5}}2$ be the Golden Ratio, and
$\phi^\sigma=\frac{1-\sqrt{5}}2$ its Galois conjugate. Let $K$ be an
imaginary quadratic number field, with discriminant $D_K\neq -4$ (to
simplify the statement in this survey), Dedekind's zeta function
$\zeta_{K}$ and ring of integers $\OOO_K$. We define as the complexity
of a quadratic irrational $\alpha$ with Galois conjugate
$\alpha^\sigma$ the quantity
$$
h(\alpha)=\frac{2}{|\alpha-\alpha^\sigma|}\;.
$$
(See \cite[\S 4.1]{ParPau12JMD} for algebraic versions and
explanations). Let $\aaa$ be a nonzero ideal in $\OOO_K$, and let
$\Ga_{0}(\aaa)$ be the congruence subgroup
$\Big\{\pm\Big(\!\begin{array}{cc}a & b \\ c & d \end{array}
\!\Big)\in {\rm PSL}_2(\OOO_K) \;:\; c\in\aaa\Big\}$. Assume (to
simplify the statement in this survey) that $\phi^\sigma$ is not in
the $\Ga_{0}(\aaa)$-orbit of $\phi$.

\bcoro [Parkkonen-Paulin \mbox{\cite[Coro.~4.7]{ParPau12JMD}, 
\cite[\S 6]{ParPau13b}}] \label{coro:orbianchi} 
There exists $\kappa>0$ such that, as $s$ tends to $+\infty$, the
cardinality of $\{\alpha\in \Ga_{0}(\aaa) \cdot \{\phi,\phi^{\sigma}\}
\!\!\!\mod \OOO_K\;:\; h(\alpha)\leq s\} $ is equal to
$$
\frac{4\pi^2\;k_\aaa\;\ln\phi} {|D_K|\;  \zeta_{K}(2)\;
  N(\aaa)\prod_{\ppp |\aaa} \big(1+\frac{1}{N(\ppp)}\big)}\; \;
s^2\,(1+\operatorname{O}(s^{-\kappa}))\;,
$$
where $k_\aaa$ is the smallest $k\in\NN-\{0\}$ such that the $2k$-th
term of Fibonacci's sequence belongs to $\aaa$, and $\ppp$ ranges over
the prime ideals in $\OOO_K$.  \cqfd \ecoro

\subsection{Counting representations of integers by binary 
Hamiltonian forms} 
\label{subsec:hamilton}

A {\em quaternion algebra} over a field $F$ is a four-dimensional
central simple algebra over $F$.  A real quaternion algebra (that is,
a quaternion algebra over $\RR$) is isomorphic either to the algebra
of real $2\times 2$ matrices over $\RR$ or to Hamilton's quaternion
algebra $\HH$ over $\RR$, with basis elements $1,i,j,k$ as a
$\RR$-vector space, with unit element $1$ and $i^2=j^2=-1$,
$ij=-ji=k$. We define the {\em conjugate} of $x=x_0+x_1i+x_2j+x_3k$ in
$\HH$ by $\overline{x}=x_0-x_1i-x_2j-x_3k$, its {\em reduced trace} by
$\tr(x)=x+\overline{x}$, and its {\em reduced norm} by $\n(x)=
x\,\overline{x}=\overline{x}\,x$.  We refer for instance to
\cite{Vigneras80} for generalities on quaternion algebras.

A {\em binary Hamiltonian form}  is a map $f:\HH \times \HH \ra \RR$
with
$$
f(u, v) = a\n(u) + \tr(\ov u\, b\, v) + c\n(v)\;,
$$
whose {\em coefficients} $a$ and $c$ are real, and $b$ lies in
$\HH$. The {\em matrix} $M(f)$ of $f$ is the Hermitian matrix
$\Big(\!\begin{array}{cc}a& b\\ \ov b& c\end{array}\!\Big)$, so that
$f(u,v)=(\!\begin{array}{cc} \!\overline u& \overline v\!
\end{array}\!)\;
\Big(\!\begin{array}{cc}a& b\\\ov b& c\end{array}\!\Big)\;
\Big(\!\begin{array}{c} \!u\!\\ \!v\!\end{array}\Big)$.  The {\em
  discriminant} of $f$ is
$$
\Delta (f) = \n(b)- ac\,,
$$ 
and $f$ is {\em indefinite} (that is, $f$ takes both positive and
negative values) if and only if $\Delta>0$.

In this section, we will describe the results in \cite{ParPau12ANT} on
the representation of integers by indefinite binary Hamiltonian forms.
The proof follows the same ideas as in the previous two sections but
the noncommutativity of the quaternions adds several new features.

In order to generalise the results of the previous two subsections to
the context of Hamiltonian forms, we have to introduce the correct
analogs of the ring of integers and of Bianchi groups for quaternion
algebras.  We say that a quaternion algebra $A$ over $\QQ$ is {\em
  definite} (or ramified over $\RR$) if the real quaternion algebra
$A\otimes_\QQ\RR$ is isomorphic to $\HH$. We fix an identification
between $A\otimes_\QQ\RR$ and $\HH$, so that $A$ is a $\QQ$-subalgebra
of $\HH$.  The {\em reduced discriminant} $D_A$ of $A$ is the product
of the primes $p\in\NN$ such that the quaternion algebra $A\otimes_\QQ
\QQ_p$ over $\QQ_p$ is a division algebra.  For example, the
$\QQ$-vector space $\HH_{\QQ}=\QQ+\QQ i+\QQ j+\QQ k$ generated by
$1,i,j,k$ in $\HH$ is Hamilton's quaternion algebra over $\QQ$.  It is
the unique (up to isomorphism) definite quaternion algebra over $\QQ$
with discriminant $D_A=2$.

A {\em $\ZZ$-lattice} $I$ in $A$ is a finitely generated $\ZZ$-module
generating $A$ as a $\QQ$-vector space.  An {\em order} in a
quaternion algebra $A$ over $\QQ$ is a unitary subring $\OOO$ of $A$
which is a $\ZZ$-lattice, and the order is {\em maximal} if it is
maximal with respect to inclusion among all orders of $A$.  The {\em
  Hurwitz order} $\OOO=\ZZ+\ZZ i+\ZZ j+ \ZZ\frac{1+i+j+k}{2}$ in
$\HH_{\QQ}$ is maximal, and it is the unique maximal order in
$\HH_{\QQ}$ up to conjugacy.

The Dieudonn\'e determinant (see \cite{Dieudonne43,Aslaksen96}) $\Det$
is the group morphism from the group $\operatorname{GL}_2(\HH)$ of
invertible $2\times 2$ matrices with coefficients in $\HH$ to
$\RR^*_+$, defined by
$$
\Det\big(\Big(\!\begin{array}{cc}a& b\\c& d\end{array}\!\Big)\big)^2\;=
\n(a\,d)+ \n(b\,c) - \tr(a\,\ov c\,d\,\ov b) =\left\{\begin{array}{cl}
\n(ad - aca^{-1}b) &{\rm if}\; a \neq 0\\
\n(cb - cac^{-1}d) &{\rm if}\; c \neq 0\\
\n(cb - db^{-1}ab) &{\rm if}\; b \neq 0\;. \end{array}\right.
$$
We will denote by $\SLH$ the group of $2\times 2$ matrices with
coefficients in $\HH$ with Dieudonn\'e determinant $1$, which equals
the group of elements of (reduced) norm $1$ in the central simple
algebra $\M_2(\HH)$ over $\RR$, see \cite[\S 9a]{Reiner75}. We
refer for instance to \cite{Kellerhals03} for more information on
$\SLH$.

The group $\SLH$ acts linearly on the left on the right $\HH$-module
$\HH\times\HH$. Let $\PP^1_r(\HH)= (\HH\times\HH-\{0\})/
\HH^\times$ be the right projective line of $\HH$, identified as usual
with the Alexandrov compactification $\HH\cup\{\infty\}$ where $[1:0]
=\infty$ and $[x:y]=xy^{-1}$ if $y\neq 0$. The projective action of
$\SLH$ on $\PP^1_r(\HH)$, induced by its linear action on $\HH\times
\HH$, is then the action by homographies on $\HH\cup\{\infty\}$
defined by
$$
\Big(\begin{array}{cc} a & b \\ c & d\end{array}\Big)\cdot z =
\left\{\begin{array}{ll} (az+b)(cz+d)^{-1} &
{\rm if}\; z\neq \infty,-c^{-1}d \\
ac^{-1} & {\rm if}\; z=\infty, c\neq 0\\
\infty & {\rm otherwise~.}\end{array}\right.
$$
The linear action on the left on $\HH\times\HH$ of the group $\SLH$
induces an action on the right on the set of binary Hermitian forms
$f$ by precomposition.

The above action of $\SLH$ on $\HH\cup\{\infty\}$ induces a faithful
left action of $\PSLH=\SLH/\{\pm \id\}$ on $\HH\cup\{\infty\}
=\partial_\infty\hcr$.  By Poincar\'e's extension procedure (see for
instance \cite[Lem.~6.6]{ParPau10GT}), this action extends to a left
action of $\SLH$ on the upper halfspace model of $\hcr$ with
coordinates $(z,r)\in\HH\times\;]0,+\infty[\,$, by
$$
 \Big(\begin{array}{cc} a & b \\ c
    & d\end{array}\Big)\cdot (z,r)= 
\Big(\;\frac{(az+b)\,\overline{(cz+d)}+a\,\overline{c}\,r^2}
{\n(cz+d)+r^2\n(c)}, \, \frac r{\n(cz+d)+r^2\n(c)}\,\Big)\;.
$$
In this way, the group $\PSLH$ is identified with the group of
orientation preserving isometries of $\hcr$.

Given an order $\OOO$ in a definite quaternion algebra over $\QQ$, a
binary Hamiltonian form $f$ is {\em integral} over $\OOO$ if its
coefficients belong to $\OOO$.  Note that such a form $f$ takes
integral values on $\OOO\times \OOO$.  The {\em Hamilton-Bianchi
  group} $\SLO= \SLH \cap\M_2(\OOO)$ preserves the set of
indefinite binary Hamiltonian forms $f$ that are integral over
$\OOO$. The stabiliser in $\SLO$ of such a form $f$ is its {\em
  group of automorphs} $\operatorname{SU}_f(\OOO)$.

The Hamilton-Bianchi group $\SLO$ is a (nonuniform) arithmetic
lattice in the connected real Lie group $\SLH$ (see for instance
\cite[p.~1104]{ParPau10GT} for details). The volume of the quotient
real hyperbolic orbifold $\SLO\bs\hcr$ has a nice expression in
terms of the discriminant $D_A$, generalising Humbert's formula.

\btheo[Emery, Parkkonen-Paulin \cite{ParPau12ANT}]
$$
\covol(\SLO)= \frac{\zeta(3)\prod_{p|D_A}(p^3-1)(p-1)\,}{11520}
\;.\;\;\;\Box
$$
\etheo

This result is proved in \cite{ParPau12ANT} using two different
methods: In the Appendix of that paper, Emery (who was the first to
prove the theorem in full generality) uses Prasad's formula and we
give a different proof using the theory of Eisenstein series for
quaternions developped in \cite{KraOse90}, following Sarnak's proof in
\cite{Sarnak83} for Bianchi groups.

With $a,b,c$ the coefficients of $f$, let
$$
\C_\infty(f)=\{[u:v]\in\PP^1_r(\HH)\;:\;f(u,v)=0\}\;\;\;{\rm and}
$$
$$
\C(f)=\{(z,r)\in\HH\times\,]0,+\infty[\;:\;f(z,1)+a\,r^2=0\}\;.
$$
In $\PP^1_r(\HH)=\HH\cup\{\infty\}$, the set $\C_\infty(f)$ is the
$3$-sphere of center $-\frac{b}{a}$ and radius $\frac{\sqrt{\Delta(f)}}
{|a|}$ if $a\neq 0$, and it is the union of $\{\infty\}$ with the real
hyperplane $\{z\in\HH\;:\; \tr(\overline{z}b)+c=0\}$ of $\HH$
otherwise.  The arithmetic group $\operatorname{SU}_f(\OOO)$ acts with
finite covolume on $\C(f)$.

The action by homographies of $\SLO$ preserves the right projective
space $\PP^1_r(\OOO) = A\cup\{\infty\}$, which is the set of fixed
points of the parabolic elements of $\SLO$ acting on
$\hcr\cup\partial_\infty \hcr$. In order to describe the orbits of
parabolic fixed points, we recall some basic definitions and facts on
ideals in a quaternion algebra, see \cite{Vigneras80}.  The {\em left
  order} $\OOO_\ell(I)$ of a $\ZZ$-lattice $I$ is $\{x\in A\;:\; xI
\subset I\}$.  A {\em left fractional ideal} of $\OOO$ is a
$\ZZ$-lattice of $A$ whose left order is $\OOO$. A {\em left ideal} of
$\OOO$ is a left fractional ideal of $\OOO$ contained in $\OOO$.  Two
 left fractional ideals $\mmm$ and $\mmm'$ of $\OOO$ are
isomorphic as left $\OOO$-modules if and only if $\mmm'=\mmm c$ for
some $c\in A^\times$. A (left) {\em ideal class} of $\OOO$ is an
equivalence class of  left fractional ideals of $\OOO$ for this
equivalence relation. We will denote by $_\OOO\!\I$ the set of ideal
classes of $\OOO$.  The {\em class number} $h_A$ of $A$ is the number
of ideal classes of a maximal order $\OOO$ of $A$. It is finite and
independent of the maximal order $\OOO$ (see for instance
\cite[p.~87-88]{Vigneras80}).

For every $(u,v)$ in $\OOO\times\OOO-\{(0,0)\}$, consider the two
left ideals of $\OOO$
$$
I_{u,v}=\OOO u+\OOO v\;,\;\;K_{u,v}=\Big\{\begin{array}{cl}
\OOO u\cap\OOO v &{\rm if~} uv\neq 0\;,\\ 
\OOO&{\rm otherwise.}\end{array}
$$
The map
$$
\SLO\bs\PP^1_r(\OOO)\ra (_\OOO\!\I\times \;_\OOO\!\I)\;,
$$
which associates, to the orbit of $[u:v]$ in $\PP^1_r(\OOO)$ under
$\SLO$, the couple of ideal classes $([I_{u,v}],[K_{u,v}])$ is a
bijection by \cite[Satz 2.1, 2.2]{KraOse90}. In particular, the number
of cusps of $\SLO$ (or the number of ends of $\SLO\bs\hcr$) is
the square of the class number $h_A$ of $A$.

The {\em norm} $\n(\mmm)$ of a left ideal $\mmm$ of $\OOO$ is the
greatest common divisor of the norms of the nonzero elements of
$\mmm$. In particular, $\n(\OOO)=1$.  The {\em norm} of a  left
fractional ideal $\mmm$ of $\OOO$ is $\frac{\n(c\mmm)}{\n(c)}$ for any
$c\in \NN-\{0\}$ such that $c\mmm \subset \OOO$.

Let $\OOO$ be a maximal order in $A$, and let $\mmm$ be a left
fractional ideal of $\OOO$, with norm $\n(\mmm)$.  For every $s>0$, we
consider the integer
$$
\psi_{f,\mmm}(s)=\card\;\;_{\operatorname{SU}_f(\OOO)}\bs
\big\{(u,v)\in\mmm\times\mmm\;:\;\n(\mmm)^{-1}|f(u,v)|\leq s,
\;\;\;\OOO u+\OOO v=\mmm\big\}\;,
$$
which is the number of nonequivalent $\mmm$-primitive representations
by $f$ of rational integers with absolute value at most $s$. 
Analogously to the cases of binary quadratic and Hermitian forms, we
have an explicit asymptotic result for this counting function.

\btheo [Parkkonen-Paulin \cite{ParPau12ANT,ParPau13b}] 
\label{theo:mainintro} There exists $\kappa>0$ such that, as $s$ tends
to $+\infty$, with $p$ ranging over positive rational primes,
$$
\psi_{f,\mmm}(s)= \frac{540\;h_A\;\covol(\operatorname{SU}_f(\OOO))}
{\pi^2\;\zeta(3)\;\Delta(f)^2\;\prod_{p|D_A}(p^3-1)(1-p^{-1})}\;\;
s^4\,(1+\operatorname{O}(s^{-\kappa}))\,.
\;\;\;\Box
$$
\etheo 

The proof of the above result again uses Corollaire 4.9 of
\cite{ParPau12JMD}, and \cite[\S 6]{ParPau13b} (see \S
\ref{sec:errorterm}) for the error term. One considers the $h_A$
different orbits of the parabolic fixed points $xy^{-1}$ of $\SLO$ for
which $\OOO x+\OOO y=\mmm$, and connects the counting functions
$$
\psi_{f,x,y}(s)=\card\;\;_{\operatorname{SU}_f(\OOO)}\bs
\big\{(u,v)\in \SLO(x,y)\;:\;
\n(\OOO x+\OOO y)^{-1}|f(u,v)|\leq s\big\}\;
$$
with the geometric counting function that counts the common
perpendiculars between a Margulis cusp neighbourhood of the cusp
corresponding to $xy^{-1}$ and the totally geodesic immersed
hypersurface corresponding to $\C(f)$.  The counting function
$\psi_{f,x,y}$ depends (besides $f$) only on the $\SLO$-orbit of
$[x:y]$ in $\PP^1_r(\OOO)$, and summing over all such orbits gives the
result. We refer to \cite{ParPau12ANT} for more details and more
general results that cover finite index subgroups of $\SLO$, and
\cite[\S 6]{ParPau13b} for the error term.

\section{Patterson, Bowen-Margulis and skinning
  measures}
\label{sec:skinning}

Let $M$ be a complete nonelementary connected Riemannian manifold of
dimension at least $2$, with pinched negative sectional curvature
$-b^2\le K\le -1$.  Let $\wt M\ra M$ be a universal Riemannian cover,
and let $\Ga$ be its covering group. Let $x_0\in \wt M$ and let
$\delta= \delta_\Ga \in\; ]0,+\infty[$ be the critical exponent of
$\Ga$.

\subsection{Patterson densities and Bowen-Margulis measures}
\label{subsec:pattbowmarg}

Let $r>0$. A family $(\mu_{x})_{x\in \wt M}$ of nonzero finite
measures on $\partial_{\infty}\wt M$ whose support is the limit set
$\Lambda\Ga$ is a {\em Patterson density of dimension} $r$ for $\Ga$
if it is $\Ga$-equivariant, that is, if it satisfies
\begin{equation}\label{eq:equivardensity}
\ga_*\mu_x=\mu_{\ga x}
\end{equation}
for all $\ga\in \Ga$ and $x\in\wt M$, and if the pairwise
Radon-Nikodym derivatives of the measures $\mu_x$ for $x\in \wt M$
exist and satisfy
\begin{equation}\label{eq:RNdensity}
\frac{d\mu_{x}}{d\mu_{y}}(\xi)=e^{-r\beta_{\xi}(x,y)}
\end{equation}
for all $x,y\in\wt M$ and  $\xi\in\partial_{\infty}\wt M$.

If {\em Poincar\'e's series}
$$
\P_\Ga(s)=\sum_{\ga\in\Ga}e^{-sd(x_0,\ga x_0)}
$$
diverges at $s=\delta$, then $\Ga$ is said to be of {\em divergence
  type}. In particular, this holds when $\wt M$ is a symmetric space
and $\Ga$ is geometrically finite by \cite{Sullivan84, CorIoz99}, see
\cite{DalOtaPei00} for many more general results. For groups of divergence
type, there exists (see for instance \cite[Coro.~1.8]{Roblin03}), up to
multiplication by a constant, one and only one Patterson density
$(\mu_x)_{x\in \wt M}$ of dimension $\delta$ for $\Ga$: For every
$x\in \wt M$, the measure $\mu_x$ is the weak-* limit of
$$
\frac{1}{\P_\Ga(s)}\;\sum_{\ga\in\Ga}e^{-s\,d(x,\ga x_0)}\Delta_{\ga x_0}
$$
as $s\to\delta$, see \cite{Patterson76, Kaimanovich90}, where
$\Delta_y$ is the unit mass Dirac measure at any point $y\in\wt M$.

Let $(\mu_{x})_{x\in \wt M}$ be a Patterson density of dimension
$\delta$ for $\Ga$. The {\em Bowen-Margulis measure} $\wt m_{\rm BM}$
for $\Ga$ on $T^1\wt M$ is defined, using Hopf's parametrisation, by
$$
d\wt m_{\rm BM}(v)=\frac{d\mu_{x_{0}}(v_{-})d\mu_{x_{0}}(v_{+})dt}
     {d_{x_{0}}(v_{-},v_{+})^{2\delta}}
=e^{-\delta(\beta_{v_{-}}(\pi(v),\,x_{0})+
\beta_{v_{+}}(\pi(v),\,x_{0}))}d\mu_{x_{0}}(v_{-})d\mu_{x_{0}}(v_{+})dt\,,  
$$
see \cite{Sullivan79, Sullivan84, Kaimanovich90}.  The Bowen-Margulis
measure is independent of the base point $x_{0}$, and its support is
(in Hopf's parametrisation) $(\Lambda\Ga\times \Lambda\Ga- \Delta)
\times\RR$, where $\Delta$ is the diagonal in $\Lambda\Ga\times\Lambda
\Ga$.  It is invariant under the geodesic flow, the antipodal map and
the action of $\Ga$, and thus it defines a measure $m_{\rm BM}$ on
$T^1M$ which is invariant under the geodesic flow of $M$ and the
antipodal map.

When the Bowen-Margulis measure $m_{\rm BM}$ is finite, the group
$\Ga$ is of divergence type (see for instance \cite[p.~19]{Roblin03}),
hence denoting the total mass of a measure $m$ by $\|m\|$, the
probability measure $\frac{m_{\rm BM}} {\|m_{\rm BM}\|}$ is then
uniquely defined, and is the unique probability measure of maximal
entropy of the geodesic flow (see \cite{OtaPei04}).  When finite, the
Bowen-Margulis measure $m_{\rm BM}$ on $T^1M$ is mixing for the
geodesic flow, under the mild assumption conjecturally always
satisfied, that the geodesic flow is topologically mixing (or that the
set of the lengths of the closed geodesics in $M$ is not contained in a
discrete subgroup of $\RR$), see \cite{Babillot02b}. This condition
holds for instance if $M$ is locally symmetric or if $M$ is compact or
if $\Ga$ contains a parabolic element, see for instance
\cite{DalBo99}.  In this review, we assume that $m_{\rm BM}$ is
finite.

\subsection{Skinning measures}

\medskip Let $\wt C$ be a nonempty closed convex subset of $\wt M$. We
define in \cite{ParPau13a} the {\em skinning measure} $\wt\sigma_{\wt C}$
of $\Ga$ on $\normal{\wt C}$, using the homeomorphism $w\mapsto w_+$ from
$\normal {\wt C}$ to $
\partial_{\infty}\wt M-\partial_{\infty}\wt C$, by
\begin{align}
d\wt\sigma_{\wt C}(w) &=  e^{-\delta\,\beta_{w(+\infty)}(\pi(w),\,x_{0})}\,d(
P^+_{\wt C})_{*}(\mu_{x_{0}}|_{\partial_{\infty}\wt 
  M-\partial_{\infty}\wt C})(w)\nonumber\\
  & =  e^{-\delta\,\beta_{w_{+}}(P_{\wt C}(w_+),\,x_{0})}\,d\mu_{x_{0}}(w_{+})\,.
\label{eq:defiskinmeas}
\end{align}
We also consider $\wt\sigma_{\wt C}$ as a measure on $T^1\wt M$ with
support contained in $\normal {\wt C}$. The skinning measure
$\wt\sigma_{\wt C}$ is independent of the base point $x_{0}$, satisfies
$\wt\sigma_{\ga \wt C} =\ga_{*}\wt\sigma_{\wt C}$ for every isometry
$\ga$ of $\wt M$ and its support is $\{w\in\normal {\wt C}:
w_{+}\in\Lambda\Gamma\}= P^+_{\wt C}(\Lambda \Gamma- \Lambda\Gamma
\cap \partial_\infty {\wt C})$.  For any $x\in\wt M$, up to
identifying the unit tangent sphere $T^1_{x}\wt M$ at $x$ with the
boundary at infinity $\partial_{\infty}\wt M$ by the map $v\mapsto
v_{+}$, we have $\wt \sigma_{\{x\}}=\mu_x$.

The skinning measure has been defined by Oh and Shah \cite[\S
1.2]{OhSha12} for the outer unit normal bundles of spheres,
horospheres and totally geodesic subspaces in real hyperbolic spaces,
see also \cite[Lem.~4.3]{HerPau10} for a closely related measure. The
terminology comes from McMullen's proof of the contraction of the
skinning map (capturing boundary information for surface subgroups of
$3$-manifold groups) introduced by Thurston to prove his
hyperbolisation theorem.

When ${\wt C}$ is a horoball, the skinning measure of ${\wt C}$ is
well known. In fact, the outer unit normal bundle $\normal {\wt C}$ of
${\wt C}$ is a leaf of the strong unstable foliation of the geodesic
flow and the skinning measure $\wt \sigma_{\wt C}$ is the conditional
measure of the Bowen-Margulis measure on this leaf, see for example
\cite{Margulis04,Roblin03}.  The skinning measure of a horoball has also
appeared as a measure on $\partial_\infty\wt M$ with the point at
infinity $\xi$ of the horoball removed in \cite{Cosentino99, Tukia94,
  Ala-Mattila11} in the constant curvature case and in \cite{HerPau04}
under the name {\em Patterson measure on} $\partial_\infty\wt M-
\{\xi\}$ in the general case.  Furthermore, using the upper halfspace
model of $\hnr$, Oh and Shah consider in \cite{OhSha12} a measure
$\omega_\Ga$ defined in $\RR^{n-1}= \partial_\infty \hnr- \{\infty\}$
by
$$
d\omega_\Ga(\xi)=e^{\delta \beta_\xi(x,\,(\xi,1))}d\mu_x(\xi).
$$
Noticing that $(\xi,1)=P_{\wt C}(\xi)$ if ${\wt C}$ is the horoball in
$\hnr$ that consists of the points whose vertical coordinate is at least
1, it follows that $\omega_\Ga$ is the image of the skinning measure
of ${\wt C}$ 
under the map $P_{\wt C}^{-1}:(\xi,1)\mapsto \xi$.

For later use in \S \ref{sec:mainstatement}, we introduce some
convenient notation. Let $w\in T^1\wt M$. When ${\wt C}=H\!B_-(w)$ is the
unstable horoball of $w$, the conditional measure of the
Bowen-Margulis measure on the strong unstable leaf $W^{\rm su}(w)$ of
$w$ is denoted by
$$
\mu_{W^{\rm su}(w)}=\wt \sigma_{H\!B_-(w)},
$$
and similarly, we denote by
$$
\mu_{W^{\rm ss}(w)}=\iota_*\wt \sigma_{H\!B_{+}(w)}
$$
the conditional measure of the Bowen-Margulis measure on the strong
stable leaf $W^{\rm ss}(w)$ of $w$. These two measures are independent
of the element $w$ of a given strong unstable leaf and given strong
stable leaf, respectively.  We also define the conditional measure of
the Bowen-Margulis measure on the stable leaf $W^{\rm s}(w)$ of $w$,
using the homeomorphism $(v',t)\mapsto v=\flow t v'$ from $W^{\rm ss}
(w) \times\RR$ to $W^{\rm s}(w)$, by
$$
d\mus{w}(v)=e^{-\delta_\Ga t}\;d\mu_{W^{\rm ss}(w)}(v')dt\;.
$$

Let ${\wt C}$ be a proper nonempty closed convex subset of $\wt M$
such that the $\Ga$-orbit of ${\wt C}$ is locally finite, and let $C$
be its image in $M$.  Since $\wt\sigma_{\wt C}$ is invariant under the
stabiliser $\Ga_{\wt C}$ of ${\wt C}$ in $\Ga$, the measure $\wt\sigma
= \sum_{\ga \in \Ga/\Ga_{\wt C}} \;\ga_*\wt\sigma_{\wt C}$ is a
$\Ga$-invariant locally finite Borel positive measure on $T^1\wt M$
(independent of the choice of representatives of elements of
$\Ga/\Ga_{\wt C}$), whose support is contained in the $\Ga$-orbit of
$\normal {\wt C}$. Hence $\wt\sigma$ induces a locally finite Borel
positive measure $\sigma_C$ on $T^1M=\Ga\bs T^1\wt M$, called the {\em
  skinning measure} of the properly immersed closed convex subset $C$,
whose support is contained in $\normal C$.

Oh and Shah proved in particular that $\|\sigma_C\|$ is finite if $\wt
M$ is geometrically finite with constant curvature $-1$ and either
${\wt C}$ is a horoball centered at a parabolic fixed point or
$\delta_\Ga>1$ and $\wt C$ is a codimension $1$ totally geodesic
submanifold. See \cite[\S 5]{OhShaCircles} for a precise, more
general statement in higher codimension. Extending this result in
variable curvature (with a different proof), we give a sharp criterion
in \cite[Theo.~9]{ParPau13a} for the finiteness of the skinning
measure, by studying its decay in the cusps of $M$. This decay is
analogous to the decay of the Bowen-Margulis measure in the cusps,
which was first studied by Sullivan \cite{Sullivan84} who called it
the fluctuating density property (see also \cite{StrVel95} and
\cite[Theo.~4.1]{HerPau04}). The criterion, as in the case of the
Bowen-Margulis measure in \cite{DalOtaPei00}, is a separation property
of critical exponents.

\subsection{Disintegration of the Bowen-Margulis measure}
\label{subsec:disintegr}

Let ${\wt C}$ be a proper nonempty closed convex subset of $\wt
M$. Define
\begin{equation}\label{eq:defiUC}
  U_{\wt C}=\{v\in T^1\wt M:\ v_+\notin\partial_\infty {\wt C}\}\;,
\end{equation}
which is a nonempty open subset of $T^1\wt M$, invariant under the
geodesic flow.

\medskip \noindent
\begin{minipage}{9.4cm} ~~~Let $f_{\wt C}:U_{\wt C}\ra \normal {\wt
    C}$ be the composition of the map from $U_{\wt C}$ onto
  $ \partial_\infty \wt M-\partial_\infty {\wt C}$ sending $v$ to
  $v_+$ and the homeomorphism $P^+_{{\wt C}}$ from $\partial_\infty
  \wt M-\partial_\infty {\wt C}$ to $\normal {\wt C}$.  The map
  $f_{\wt C}$ is a continuous fibration, invariant under the geodesic
  flow. The fiber of $f_{\wt C}$ above $w\in \normal {\wt C}$ is
  exactly the stable leaf $W^{\rm s}(w)=\{v\in T^1\wt M\;:\;v_+
  =w_+\}$.  See \cite{ParPau13a,ParPau13b} for further properties of
  $f_{\wt C}$, including the fact that $f_{\wt C}$ is a H\"older
  fibration when the sectional curvature of $M$ has bounded
  derivatives.
\end{minipage}
\begin{minipage}{5.5cm}
\begin{center}
\begin{picture}(0,0)%
\includegraphics{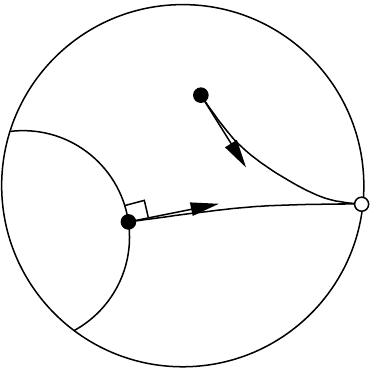}%
\end{picture}%
\setlength{\unitlength}{3812sp}%
\begingroup\makeatletter\ifx\SetFigFont\undefined%
\gdef\SetFigFont#1#2#3#4#5{%
  \reset@font\fontsize{#1}{#2pt}%
  \fontfamily{#3}\fontseries{#4}\fontshape{#5}%
  \selectfont}%
\fi\endgroup%
\begin{picture}(1860,1814)(1613,-1688)
\put(2397,-1094){\makebox(0,0)[lb]{\smash{{\SetFigFont{11}{13.2}{\rmdefault}{\mddefault}{\updefault}{\color[rgb]{0,0,0}$f_{\wt C}(v)$}%
}}}}
\put(1846,-781){\makebox(0,0)[lb]{\smash{{\SetFigFont{11}{13.2}{\rmdefault}{\mddefault}{\updefault}{\color[rgb]{0,0,0}$\wt C$}%
}}}}
\put(2341,-1321){\makebox(0,0)[lb]{\smash{{\SetFigFont{11}{13.2}{\rmdefault}{\mddefault}{\updefault}{\color[rgb]{0,0,0}$=P^+_{\wt C}(v_+)$}%
}}}}
\put(3458,-935){\makebox(0,0)[lb]{\smash{{\SetFigFont{11}{13.2}{\rmdefault}{\mddefault}{\updefault}{\color[rgb]{0,0,0}$v_+$}%
}}}}
\put(2713,-455){\makebox(0,0)[lb]{\smash{{\SetFigFont{11}{13.2}{\rmdefault}{\mddefault}{\updefault}{\color[rgb]{0,0,0}$v$}%
}}}}
\end{picture}%

\end{center}
\end{minipage}

\bigskip The following disintegration result of the Bowen-Margulis
measure over the skinning measure of $\wt C$ is the crucial tool for the
proof in \cite{ParPau13b} of our general counting result, see \S
\ref{sec:mainstatement}.

\bprop[Parkkonen-Paulin \cite{ParPau13a}] \label{prop:disintegration}
Let $\wt C$ be a proper nonempty closed convex subset of $\wt M$. The
restriction to $U_{\wt C}$ of the Bowen-Margulis measure $\wt m_{\rm
  BM}$ disintegrates by the fibration $f_{\wt C}:U_{\wt C}\ra\normal
{\wt C}$, over the skinning measure $\wt\sigma_{{\wt C}}$ of ${\wt
  C}$, with conditional measure 
$e^{\delta\;  \beta_{w_+} (\pi(w),\,\pi(v))}\;d\mus{w}(v)$ 
on the fiber $f_{\wt C}^{-1}(w)
=W^s(w)$ of $w\in \normal {\wt C}$: 
$$
d\wt m_{\rm BM}(v)=\int_{w\in \normal {\wt C}}
\;e^{\delta\; \beta_{w_+}(\pi(w),\,\pi(v))}\;
d\mus{w}(v)\;d\wt\sigma_{{\wt C}}(w)\;.\;\;\;\Box
$$
\eprop

\section{Finite volume hyperbolic manifolds}
\label{sec:finitevolhypman}

In this section, we consider the special case when $\wt M=\hnr$, $\Ga$
is a discrete group of isometries of $\wt M$ and $M=\Ga\bs\wt M$ has
finite volume, and we relate the measures defined in \S
\ref{sec:skinning} with more classical measures. For every $p \in\NN$,
we denote by $\lambda_{p}$ the standard Lebesgue measure of $\RR^{p}$.

Under the assumptions of this section, there exists a unique Patterson
density $(\mu_x)_{x\in\hnr}$ of dimension $n-1$ for $\Ga$ normalised
to have total mass $\Vol(\SSS^{n-1})$ for every $x\in\hnr$, which we
call the spherical density and which we now describe.

In the unit ball model of $\hnr$ with origin $0$, the measure $\mu_0$
of the {\em spherical density} $(\mu_x)_{x\in\hnr}$ is the Lebesgue
measure of $\SSS^{n-1}=\partial_\infty\hnr$ and (see for instance
\cite[p.~273]{BriHae99})
$$
\frac{d\mu_x}{d\mu_0}(\xi)=e^{-(n-1)\beta_{\xi}(x,0)}=
\left(\frac{1-\|x\|^2}{\|x-\xi\|^2}\right)^{n-1}\;. 
$$
In the upper halfspace model with point at infinity $\infty$, using
the standard inversion mapping the ball model to the upper halfspace
model, the {\em spherical density} $(\mu_x)_{x\in\hnr}$ has the
expression, for every $\xi\in
\RR^{n-1}=\partial_\infty\hnr-\{\infty\}$,
\begin{equation}\label{eq:spericaldensityupper}
d\mu_x(\xi)= \big(\frac {2
  x_n}{\|x-\xi\|^2}\big)^{n-1}\;d\lambda_{n-1}(\xi)\;,
\end{equation}
where  $x_n$ is the vertical coordinate of any $x\in\hnr$.

In the unit ball model of $\hnr$, the visual distance $d_{0}$ seen
from the origin $0$ (see \S \ref{subsec:confstruc}) coincides
with half the chordal distance (see for example \cite{Bourdon95}). In
the upper halfspace model, an easy computation shows that the Busemann
cocycle of $\hnr$ is
\begin{equation}\label{eq:busemancocyclupper}
  \beta_\xi(x,y)=
  \ln(\frac {y_{n}}{x_{n}}\frac{\|x-\xi\|^2}{\|y-\xi\|^2})
\end{equation}
for all $x,y\in\hnr$ and all $\xi\in\RR^{n-1}$. By Equation
\eqref{eq:defdistvistroi}, for any base point $x\in\hnr$ and all
$\xi,\eta\in\RR^{n-1}$, using the point $u= (\frac{\xi+\eta}{2},
\frac{\|\xi-\eta\|}{2})$ as a chosen point on the geodesic line with
endpoints $\xi,\eta$, we get an expression for the visual distance
seen from $x$:
$$
d_{x}(\xi,\eta)=
\frac{x_{n}\|\xi-\eta\|}{\|x-\xi\|\|x-\eta\|}\;.
$$
Thus, in the upper halfspace model, for any $v\in T^1\hnr$ such that
$v_\pm\neq \infty$, we have
\begin{equation}\label{eq:bowenmargupper}
d\wt m_{\rm BM}(v)=\frac{2^{2(n-1)}d\lambda_{n-1}(v_{-})\;
  d\lambda_{n-1}(v_{+})\;dt}{\|v_+-v_-\|^{2(n-1)}}\;,
\end{equation}
where $t$ is the signed distance from the closest point to $\infty$ on
the geodesic line $]v_-,v_+[$ to $\pi(v)$.

It is known that the Liouville measure, normalised to be a probability
measure, is the probability measure of maximal entropy for the
geodesic flow in constant curvature and finite volume. Thus, the
Bowen-Margulis measure coincides (up to a positive multiplicative
constant) with the Liouville measure. We now determine the
proportionality constant.  

\bprop Let $M$ be a finite volume complete hyperbolic manifold of
dimension $n\geq 2$, $d\Vol_{T^1M}$ its Liouville measure, and $dm_{\rm BM}$
its Bowen-Margulis measure, constructed  using the
spherical Patterson density. Then
$$
m_{\rm BM}=2^{n-1}\,\Vol_{T^1M}\;.
$$
\eprop

In particular, 
\begin{equation}\label{eq:relatbowmargliouv}
\|m_{\rm BM}\|=2^{n-1}\Vol(\SSS^{n-1})\Vol(M)\;.
\end{equation}

\medskip
\dem We use the upper halfspace model 
$$
\hnr=\{x=(\ov x,x_n)\in\RR^{n-1}\times \RR \;:\;x_n>0\}\;.
$$
We parametrise the unit tangent sphere at any point $x\in \hnr$ by the
positive endpoint $v_+\in\RR^{n-1}\cup\{\infty\}$ of a unit tangent
vector $v\in T^1_x\hnr$. This gives a parametrisation of the vectors
$v\in T^1\hnr$ by the pairs $(x,v_+)\in \hnr\times (\RR^{n-1} \cup
\{\infty\})$.  Recall that the Liouville measure disintegrates as
$$
d\Vol_{T^1\hnr}(v)=\int_{x\in \hnr} d\Vol_{T^1_x\hnr}(v)\;d\Vol_{\hnr}(x)\;.
$$ 
Hence in the full-measure subset where $v_+\ne\infty$, the Liouville
measure may be written
\begin{equation}\label{eq:liouville}
d\Vol_{T^1\hnr}(v)=\frac{(2x_n)^{n-1}\; 
d\lambda_{n-1}(v_+)}{\|x-v_+\|^{2(n-1)}}\;\frac{d\lambda_n(x)}{x_n^n} 
=\frac{2^{n-1}d\lambda_{n-1}(\ov x)\,d\lambda_{n-1}(v_+)dx_n}
{\|x-v_+\|^{2(n-1)}\;x_n}\;.
\end{equation}

In order to relate the formulas \eqref{eq:bowenmargupper} and
\eqref{eq:liouville}, let us give the expression of the coordinates
$(\ov x,x_n,v_+)$ in terms of the coordinates $(v_-,v_+,t)$.

\medskip\noindent
\begin{minipage}{8.4cm}~~~ Let $\alpha$ be the angle between the
  segments $[\frac{v_-+v_+}{2},v_+]$ and $[\frac{v_-+v_+}{2},x]$. Let
  $\rho$ be the algebraic distance from $\frac{v_-+v_+}{2}$ to $\ov x$
  on the line through $v_-$ and $v_+$ oriented from $v_-$ to $v_+$. We
  have
$$
\ov x= \frac{v_-+v_+}{2}+\rho\frac{v_+-v_-}{\|v_+-v_-\|} 
$$
and by a formula of \cite[p.~147]{Beardon83},
$$
\sinh t= \frac{1}{\tan \alpha} =\frac{\rho}{x_n}\;.
$$
\end{minipage}
\begin{minipage}{6.5cm}
\begin{center}
\begin{picture}(0,0)%
\includegraphics{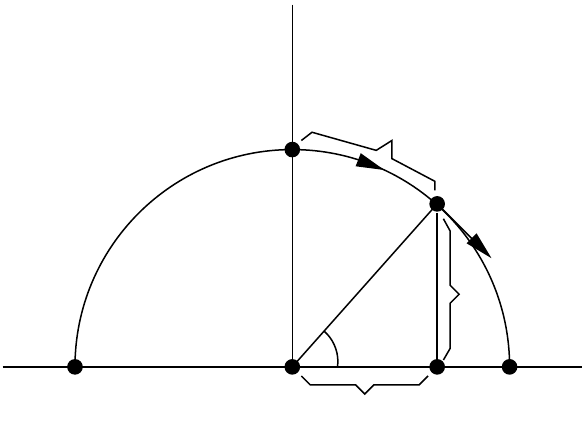}%
\end{picture}%
\setlength{\unitlength}{3812sp}%
\begingroup\makeatletter\ifx\SetFigFont\undefined%
\gdef\SetFigFont#1#2#3#4#5{%
  \reset@font\fontsize{#1}{#2pt}%
  \fontfamily{#3}\fontseries{#4}\fontshape{#5}%
  \selectfont}%
\fi\endgroup%
\begin{picture}(2904,2146)(1339,-3905)
\put(3286,-2401){\makebox(0,0)[lb]{\smash{{\SetFigFont{11}{13.2}{\rmdefault}{\mddefault}{\updefault}{\color[rgb]{0,0,0}$t$}%
}}}}
\put(3632,-3293){\makebox(0,0)[lb]{\smash{{\SetFigFont{11}{13.2}{\rmdefault}{\mddefault}{\updefault}{\color[rgb]{0,0,0}$x_n$}%
}}}}
\put(3826,-3751){\makebox(0,0)[lb]{\smash{{\SetFigFont{11}{13.2}{\rmdefault}{\mddefault}{\updefault}{\color[rgb]{0,0,0}$v_+$}%
}}}}
\put(3039,-3506){\makebox(0,0)[lb]{\smash{{\SetFigFont{11}{13.2}{\rmdefault}{\mddefault}{\updefault}{\color[rgb]{0,0,0}$\alpha$}%
}}}}
\put(1666,-3766){\makebox(0,0)[lb]{\smash{{\SetFigFont{11}{13.2}{\rmdefault}{\mddefault}{\updefault}{\color[rgb]{0,0,0}$v_-$}%
}}}}
\put(3578,-2743){\makebox(0,0)[lb]{\smash{{\SetFigFont{11}{13.2}{\rmdefault}{\mddefault}{\updefault}{\color[rgb]{0,0,0}$x$}%
}}}}
\put(3818,-3035){\makebox(0,0)[lb]{\smash{{\SetFigFont{11}{13.2}{\rmdefault}{\mddefault}{\updefault}{\color[rgb]{0,0,0}$v$}%
}}}}
\put(3151,-3841){\makebox(0,0)[lb]{\smash{{\SetFigFont{11}{13.2}{\rmdefault}{\mddefault}{\updefault}{\color[rgb]{0,0,0}$\rho$}%
}}}}
\put(2300,-3784){\makebox(0,0)[lb]{\smash{{\SetFigFont{11}{13.2}{\rmdefault}{\mddefault}{\updefault}{\color[rgb]{0,0,0}$\frac{v_-+v_+}{2}$}%
}}}}
\put(3552,-3748){\makebox(0,0)[lb]{\smash{{\SetFigFont{11}{13.2}{\rmdefault}{\mddefault}{\updefault}{\color[rgb]{0,0,0}$\ov x$}%
}}}}
\end{picture}%

\end{center}
\end{minipage}

\medskip 
Since $\rho^2+x_n^2= \|\frac{v_+-v_-}{2}\|^2$, we hence have
$$ 
x_n=\frac{\|v_+-v_-\|}{2\cosh t}
\;\;\;{\rm and}\;\;\;
\ov x= \frac{v_-+v_+}{2}+\frac{v_+-v_-}{2}\;\tanh t\;.
$$
Writing $\ov x=(\ov x^1,\dots,\ov x^{n-1})$ and
$v_\pm=(v_\pm^1,\dots,v_\pm^{n-1})$ and differentiating the above
equations with $v_+$ constant, we have, for $i=1\dots, n-1$,
$$
dx_n=-\frac{\sinh t}{2\cosh^2 t}\;\|v_+-v_-\|\;dt-
\frac{1}{2\cosh t}\sum_{j=1}^{n-1}
\frac{v_+^j-v_-^j}{\|v_+-v_-\|}\;dv_-^j
$$
and
$$
d\ov x^i=\frac{1-\tanh t}{2}\;dv_-^i+
\frac{v_+^i-v_-^i}{2\cosh^2 t}\;dt\;.
$$
Therefore an easy computation, using the facts that $\ov
x-v_+=\frac{1-\tanh t}{2}(v_--v_+)$ and $\|x-v_+\|^2= \|\ov
x-v_+\|^2+x_n^2=\frac{1-\tanh t}{2}\|v_+-v_-\|^2$, shows that
%\begin{align*}
%d\ov x^1\wedge\dots \wedge d\ov x^{n-1}&=
%\big(\frac{1-\tanh t}{2}\big)^{n-1}
%dv_-^1\wedge\dots \wedge dv_-^{n-1}\\+&
%\sum_{i=1}^{n-1}(-1)^{n-i-1}\big(\frac{1-\tanh t}{2}\big)^{n-2}
%\frac{v_+^i-v_-^i}{2\cosh^2 t} dv_-^1\wedge\dots
%\wedge \wh{dv_-^{i}}\wedge\dots \wedge dv_-^{n-1}\wedge dt
%\end{align*}
\begin{align*}
d\ov x^1\wedge\dots \wedge d\ov x^{n-1}\wedge dx_n&=
\frac{\|v_+-v_-\|}{2\cosh t}
\Big(\frac{1-\tanh t}{2}\Big)^{n-1}\;
dv_-^1\wedge\dots \wedge dv_-^{n-1}\wedge dt\\ &=
x_n\Big(\frac{\|x-v_+\|^2}{\|v_+-v_-\|^2}\Big)^{n-1}\;
dv_-^1\wedge\dots \wedge dv_-^{n-1}\wedge dt\;.
\end{align*}
The result then follows from the formulas \eqref{eq:bowenmargupper} and
\eqref{eq:liouville}. \cqfd

\medskip Let now $C$ be either a Margulis cusp neighbourhood in $M$ or
a totally geodesic immersed submanifold of $M$ with finite volume, and
let us relate the skinning measure of $C$ to the usual Riemannian
measure on the outer unit normal bundle of $C$. Note that the
Riemannian measure $\Vol_{\normal C}$ disintegrates with respect to
the base point fibration $\normal C\ra \partial C$ over the Riemannian
measure of $\partial C$, with measure on the fiber of $x\in\partial C$
the spherical measure on the outer unit normal vectors to $C$ at $x$:
\begin{equation}\label{eq:desinteRiemnormal}
d\Vol_{\normal C}(v)=
\int_{x\in \partial C} d\Vol_{\normal C\cap T^1_xM}(v)\;d\Vol_{\partial C}(x)\;.
\end{equation}

Homogeneity considerations show that the skinning measure $\sigma_C$
coincides up to a multiplicative constant with the Riemannian measure
$\Vol_{\normal C}$.  We now compute the constant.

\bprop \label{prop:calcskinning} Let $M$ be a finite volume complete
hyperbolic manifold of dimension $n\geq 2$. We use the spherical
Patterson density to define the skinning measures.

\noindent (1) If $C$ is a Margulis cusp neighbourhood, then
$$
\sigma_C= 2^{n-1}\Vol_{\normal C}\;.
$$

\noindent (2) If $C$ is a finite volume totally geodesic properly immersed
submanifold of $M$, then
$$
\sigma_C= \Vol_{\normal C}\;.
$$
\eprop

In particular, if $C$ is a Margulis cusp neighbourhood of $M$, then
(see for instance \cite[p.~473]{Hersonsky93} for the last equality)
$$
\|\sigma_{C}\|=2^{n-1}\Vol(\normal C)=2^{n-1}\Vol(\partial C)
=2^{n-1}(n-1)\Vol(C)\;,
$$
and if $C$ is a finite volume totally geodesic properly immersed
submanifold of dimension $k\in\{1,\dots,n-1\}$ of $M$, then
$$
\|\sigma_C\|=\Vol(\SSS^{n-k-1})\Vol(C)\;.
$$

\medskip \dem (1) Consider the horoball $\wt C$ in the upper halfspace
model of $\hnr$ that consists of the points whose vertical
coordinate is at least $1$. Fix a base point $x_{0}=(0,1)\in \RR^{n-1}
\times\;]0,+\infty[\,$.  Note that the closest point to $\xi\in
\RR^{n-1}$ in $\wt C$ is $P_{\wt C}(\xi)= (\xi,1)\in \RR^{n-1} \times
\;]0,+\infty[\,$. Using the definition of the skinning measure for the
first equality and the formulas \eqref{eq:spericaldensityupper} and
\eqref{eq:busemancocyclupper} for the second one, we hence have
\begin{align*}
 d\wt \sigma_{\wt C}(w) &=  
e^{-(n-1)\beta_{w_{+}}(P_{\wt C}(w_+),\,x_{0})}\,d\mu_{x_{0}}(w_{+}) \\ & = 
 (\|x_{0}-w_{+}\|^2)^{n-1} (\frac{2}{\|x_0-w_{+}\|^2}\big)^{n-1}
\;d\lambda_{n-1}(w_{+})=2^{n-1}d\lambda_{n-1}(w_{+})\,. 
\end{align*}
Since $\partial \wt C=\{(\ov x,1)\;:\;\ov x\in\RR^{n-1}\}$ is a
codimension one submanifold of $\hnr$, whose induced Riemannian metric
is isometric to the Euclidean metric on $\RR^{n-1}$ by the map $(\ov
x,1)\mapsto \ov x$, the result follows.

\medskip
(2) Let $1\leq k\leq n-1$. In the upper halfspace model of $\hnr$ with
base point $x_0=(0,\dots,0,1)$, consider the $k$-dimensional totally
geodesic subspace,
$$
\wt C=\{x=(x_1,\dots, x_n)\in\hnr: x_{1}=\cdots=x_{n-k}=0\}\;,
$$  
which is isometric to $\hkr$ and has Riemannian volume $d\Vol_{\wt C}=
\frac{d\lambda_{k-1}(x_{n-k+1},\dots,x_{n-1})d\lambda_{1}(x_n)}{x_n^k}$.

\smallskip \noindent
\begin{minipage}{7.4cm} ~~~ For any $\xi=(\xi^{1},\xi^{2})\in \RR^{n-k}
  \times\RR^{k-1}=\RR^{n-1} =\partial_\infty \hnr -\{\infty\}$, the
  closest point to $\xi$ in $\wt C$ is $P_{\wt C} (\xi)= (0,\xi^{2},
  \|\xi^{1}\|) \in\RR^{n-k} \times\RR^{k-1} \times \;]0,+\infty[\; =
  \hnr$. Note that $\pi(w)_n=\|w^1_+\|$ and $\|\pi(w) -w_+\|^2 =
  2\,\pi(w)_n \|w^1_+\|$ for every $w\in \normal{\wt C}$.  Recall that
  the map $w\mapsto w_+$ from $\normal {\wt C}$ to $\partial_{\infty}
  \hnr - \partial_{\infty}\wt C= \RR^{n-1} - (\RR^{k-1} \times \{0\})$
  is a homeomorphism. 
\end{minipage}
\begin{minipage}{7.5cm}
\begin{center}
\begin{picture}(0,0)%
\includegraphics{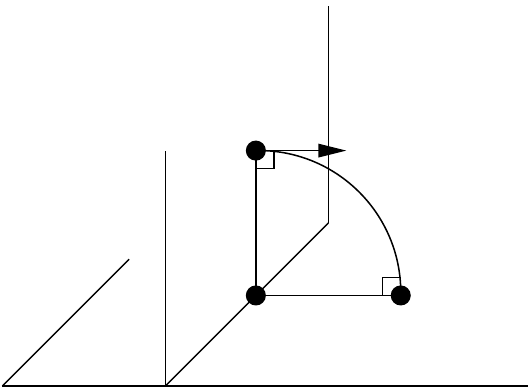}%
\end{picture}%
\setlength{\unitlength}{3812sp}%
\begingroup\makeatletter\ifx\SetFigFont\undefined%
\gdef\SetFigFont#1#2#3#4#5{%
  \reset@font\fontsize{#1}{#2pt}%
  \fontfamily{#3}\fontseries{#4}\fontshape{#5}%
  \selectfont}%
\fi\endgroup%
\begin{picture}(2634,1916)(979,-1335)
\put(2701,-196){\makebox(0,0)[lb]{\smash{{\SetFigFont{11}{13.2}{\rmdefault}{\mddefault}{\updefault}{\color[rgb]{0,0,0}$w$}%
}}}}
\put(2127,-39){\makebox(0,0)[lb]{\smash{{\SetFigFont{11}{13.2}{\rmdefault}{\mddefault}{\updefault}{\color[rgb]{0,0,0}$\pi(w)$}%
}}}}
\put(2854,-1073){\makebox(0,0)[lb]{\smash{{\SetFigFont{11}{13.2}{\rmdefault}{\mddefault}{\updefault}{\color[rgb]{0,0,0}$w_+=(w^1_+,w^2_+)$}%
}}}}
\put(2150,-1081){\makebox(0,0)[lb]{\smash{{\SetFigFont{11}{13.2}{\rmdefault}{\mddefault}{\updefault}{\color[rgb]{0,0,0}$(0,w^2_+)$}%
}}}}
\put(1843,-770){\makebox(0,0)[lb]{\smash{{\SetFigFont{11}{13.2}{\rmdefault}{\mddefault}{\updefault}{\color[rgb]{0,0,0}$\wt  C$}%
}}}}
\put(1211,-1266){\makebox(0,0)[lb]{\smash{{\SetFigFont{11}{13.2}{\rmdefault}{\mddefault}{\updefault}{\color[rgb]{0,0,0}$\RR^{n-1}$}%
}}}}
\end{picture}%

\end{center}
\end{minipage}

\medskip Using the definition of the skinning measure for the first
equality and the formulas \eqref{eq:spericaldensityupper} and
\eqref{eq:busemancocyclupper} for the second one, we hence get
\begin{align*}
d\wt\sigma_{\wt C}(w) &=
e^{-(n-1)\beta_{w_{+}}(P_{\wt C}(w_+),\,x_{0})}\,d\mu_{x_{0}}(w_{+})\\ &=
\Big(\frac{\pi(w)_n}{1}\frac{\|x_0-w_+\|^2}{\|\pi(w)-w_+\|^2}\Big)^{n-1}
\Big(\frac{2}{\|x_0-w_+\|^2}\Big)^{n-1}\;d\lambda_{n-1}(w_+)
=\frac{d\lambda_{n-1}(w_+)}{\|w^1_+\|^{n-1}}\,.
\end{align*}
On the other hand, by Equation \eqref{eq:desinteRiemnormal}, we have
$$
d\Vol_{\normal{\wt C}}(w)= 
d\Vol_{\SSS^{n-k-1}}\big(\frac{w^1_+}{\|w^1_+\|}\big)\;
\frac{d\lambda_{k-1}(w_+^2)\;d\lambda_1(\|w^1_+\|)}{\|w^1_+\|^k}\;.
$$
Using spherical coordinates on the first factor of $\RR^{n-1}=
\RR^{n-k} \times \RR^{k-1}$, we have
$$
d\lambda_{n-1}(w_+)=
\|w^1_+\|^{n-k-1} \;d\Vol_{\SSS^{n-k-1}}\big(\frac{w^1_+}{\|w^1_+\|}\big)\;
d\lambda_1(\|w^1_+\|)\; d\lambda_{k-1}(w_+^2)\;.
$$
Hence $\wt\sigma_{\wt C} =\Vol_{\normal{\wt C}}$, and the result
follows by taking quotients.  \cqfd

\section{The main counting result of common 
perpendiculars}
\label{sec:mainstatement}

Let $M$ be a nonelementary complete connected Riemannian manifold with
dimension at least $2$ and pinched sectional curvature at most
$-1$. Let $\wt M\ra M$ be a universal Riemannian cover of $M$, with
covering group $\Ga$. Let $\delta$ be the critical exponent of $\Ga$.
We assume that the Bowen-Margulis measure $m_{\rm BM}$ of $M$ is
finite and mixing for the geodesic flow.

\btheo[Parkkonen-Paulin  \cite{ParPau13b}]\label{theo:mainsurvey} 
Let $C_-$ and $C_+$ be two properly immersed clo\-sed convex subsets of
$M$. Assume that their skinning measures $\sigma_{C_-}$ and
$\sigma_{C_+}$ are finite and nonzero. Then, as $s\to+\infty$,
$$
\N_{C_-,\,C_+}(s)\sim
\frac{\|\sigma_{C_-}\|\,\|\sigma_{C_+}\|}
{\delta\,\|m_{\rm BM}\|}\;e^{\delta\, s}\,.
$$
\etheo

\noindent As in Herrmann's result (see Equation \eqref{eq:herrmann2} in
\S \ref{sec:pointandgeodesic}) or Oh-Shah's result (see the end
of \S \ref{sec:horoballgeodesic}), the endpoints of the common
perpendiculars are evenly distributed simultaneously on $C_-$ and on
$C_+$, in the following sense.

\btheo[Parkkonen-Paulin  \cite{ParPau13b}]
\label{theo:mainsurveyeffective} Let $C_-$ and $C_+$ be two properly
immersed closed convex subsets of $M$. Let $\Omega^-$ and $\Omega^+$
be relatively compact subsets of $\partial_+^1C_-$ and
$\partial_+^1C_+$, respectively. Assume that
$\sigma_{C_-}(\Omega^-)\neq 0$, $\sigma_{C_+}(\Omega^+)\neq 0$ and
$\sigma_{C_-}(\partial \Omega^-)=\sigma_{C_+}(\partial \Omega^+)=
0$. Then, as $s\to+\infty$, the number $\N_{\Omega^-,\,\Omega^+}(s)$
of common perpendiculars of $C_-$ and $C_+$, with lengths at most $s$,
and with initial vector in $\Omega^-$ and terminal vector in
$\iota\,\Omega^+$, satisfies
$$
\N_{\Omega^-,\,\Omega^+}(s)\sim
\frac{\sigma_{C_-}(\Omega^-)\;\sigma_{C_+}(\Omega^+)}
{\delta\,\|m_{\rm BM}\|}\;e^{\delta\, s}\,.
$$
\etheo

When $C_-=\{x\}$, $C_+=\{y\}$, are singletons, with $\wt x,\wt y$
lifts of $x,y$ to $\wt M$, we recover Roblin's result in
\cite{Roblin03} that
$$
\N_{C_-,\,C_+}(s)\sim
\frac{\|\mu_{\wt x}\|\,\|\mu_{\wt x}\|} {\delta\,\|m_{\rm BM}\|}\;e^{\delta\, s}
=\frac{\|\sigma_{\{x\}}\|\,\|\sigma_{\{y\}}\|}
{\delta\,\|m_{\rm BM}\|}\;e^{\delta\, s}\,.
$$

Let us give a brief sketch of proof of these results, which uses
directly the mixing property of the geodesic flow (and avoids the
equidistribution step in Margulis's scheme of proof). This will, in
particular, allow us in \S \ref{sec:errorterm} to give estimates on
the error terms in the presence of exponential decay of
correlations. We refer to \cite{ParPau13b} for complete proofs, and we
only give here a reading guide, the actual proofs require a much more
technical approach.

By definition, $C_-$ and $C_+$ are the images in $M$ of two proper
nonempty closed convex subsets $\wt C_-$ and $\wt C_+$ in $\wt M$,
whose $\Ga$-orbits are locally finite.

We introduce dynamical neighbourhoods of $\partial_+^1C_-$ and
$\partial_+^1C_+$, and we define bump functions supported in them, to
which we will apply the mixing property. We fix $\eta>0$ small enough
and $R>0$ big enough.

For every $w\in T^1\wt M$, let $V_{w,R}$ be the ball of center $w$ and
radius $R$ for Hamenst\"adt's distance $d_{W^{\rm ss}(w)}$ on the strong
stable leaf $W^{\rm ss}(w)$ of $w$ (see \S \ref{subsec:stabunstab}).  For
every proper nonempty closed convex subset $\wt D$ in $\wt M$ whose
$\Ga$-orbit is locally finite, let $\V_{\eta,R}(\wt D)$ be the union
for all $w\in \normal{\wt D}$ and $s\in\;]-\eta,\eta[$ of the sets
$g^sV_{w,R}$.  These dynamical neighbourhoods $\V_{\eta,R}(\wt D)$ of
$\normal{\wt D}$ are natural under isometries, hence, with $D$ the
image of $\wt D$ in $M$, they allow to define nice neighbourhoods
$\V_{\eta,R}(D)$ of $\partial_+^1D$, that scale nicely under the
geodesic flow: $g^t\V_{\eta,R}(\wt D)= \V_{\eta,e^{-t}R}(\N_t\wt D)$
for every $t\geq 0$.

Let $h_{\eta,\,R}: T^1\wt M\ra[0,+\infty]$ be the measurable
$\Ga$-invariant map defined by $w\mapsto\frac{1}{2\eta\;\mu_{W^{\rm
      ss}(w)} (V_{w,\,R})}$.  The constant $R>0$ is chosen big enough
so that the above denominator is nonzero if $w\in\normal{\wt C_{\pm}}$
(see \cite[Lem.~7]{ParPau13b}).  We denote by $\mathbbm{1}_A$ the
characteristic function of a subset $A$.  Let $\wt\phi_{\eta,\wt D} :
T^1\wt M\to[0,+\infty]$ be the map defined by (using the convention
$\infty\times 0=0$)
$$
\wt\phi_{\eta,\wt D}(v)=\sum_{\ga \in\Ga/\Ga_{\wt D}}
h_{\eta,\,R}\circ f_{\ga\wt D}(v)\;\mathbbm{1}_{\V_{\eta,\,R}(\ga\wt D)}(v)\;,
$$  
where $\Ga_{\wt D}$ is the stabiliser of $\wt D$ in $\Ga$ and
$h_{\eta,\,R}\circ f_{\ga\wt D}(v)\;\mathbbm{1}_{\V_{\eta,\,R}(\ga\wt
  D)}(v)=0$ if $v\notin U_{\ga \wt D}$, since $\V_{\eta,\,R}(\ga\wt
D)\subset U_{\ga \wt D}$. The function $\wt\phi_{\eta,\wt D}$ is
invariant under $\Ga$, hence defines by taking the quotient by $\Ga$ a
test function $\phi_{\eta,D}:T^1M \to[0,+\infty]$. Now define
$\phi_\eta^-= \phi_{\eta,C_-}$ and $\phi_\eta^+ =\phi_{\eta,C_+}\circ
\iota$. The invariance of the Bowen-Margulis measure by the antipodal
map and the disintegration result of Proposition
\ref{prop:disintegration} allow to prove (see
\cite[Lem.~14]{ParPau13b}) that
\begin{equation}\label{eq:intphieta}
\int_{T^1M}\phi^\pm_\eta\;dm_{\rm BM}=\|\sigma_{C_\pm}\|\;,
\end{equation}
and that $\phi^\pm_\eta\;dm_{\rm BM}\;\stackrel{*}{\longrightarrow}\;
\sigma_{C_\pm}$ as $\eta$ goes to $0$.

The main trick in the proof is to estimate in two ways the integral
$$
\I_\eta(t)=\int_{T^1M}\phi^-_\eta\circ\flow{-t/2}
\;\phi^+_\eta\circ\flow{t/2}\,dm_{\rm BM}\;.
$$
On one hand, by Equation \eqref{eq:intphieta} and the mixing property
of the geodesic flow, the integral $\I_\eta(t)$ converges, for every
fixed $\eta>0$, to $\frac{\|\sigma_{C_-}\|\;\|\sigma_{C_-}\|}{\|m_{\rm
    BM}\|}$ as $t\ra+\infty$.

On the other hand, a vector $v\in T^1M$, with a fixed lift $\wt v$ to
$T^1\wt M$, belongs to the support of $\phi^-_\eta\circ\flow{-t/2}
\;\phi^+_\eta\circ\flow{t/2}$ if and only if $\flow{-t/2}v$ belongs to
the support of $\phi^-_\eta$ and $\flow{t/2}v$ belongs to the support
of $\phi^+_\eta$, that is, if and only if there exist $\ga^\pm\in
\Ga$, $s^\pm\in\;]-\eta,\eta[$, $w^\pm\in \ga^\pm \partial^1_+ \wt
C_\pm$ and $v^\pm\in V_{w^\pm,R}$ such that $\wt v= g^{\frac{t}{2} +
  s^-} v^- =g^{-\frac{t}{2}-s^+}\iota v^+$. For every $\epsilon>0$, by
the properties of negative curvature, this implies, if $\eta$ is small
enough, and uniformly in $t$ big enough, that $\pi(\wt v)$ is not far
from the midpoint of a common perpendicular arc between $\ga^- \wt
C_-$ and $\ga^+ \wt C_+$, of length close to $t$, and that
$g^{t/2}\ga^- \partial^1_+ \wt C_-$ is close to a piece of strong
unstable leaf at $\wt v$, and $g^{-t/2}\ga^+ \iota\, \partial^1_+ \wt
C_+$ is close to a piece of strong stable leaf at $\wt v$ (see
\cite[Lem.~7]{ParPau13b}).  Furthermore, each such midpoint
contributes to the integral $\I_\eta(t)$ by an amount which is, as
$\eta$ is small and uniformly in $t$ big enough, almost
$\frac{e^{-\delta t}}{2\eta}$. By a Cesaro type of argument, the
results follows, by integrating $e^{\delta t}$.

To pass from Theorem \ref{theo:mainsurvey} to Theorem
\ref{theo:mainsurveyeffective}, we replace $\partial_+^1C_-$ and
$\partial_+^1C_+$ by $\Omega^-$ and $\Omega^+$, the endpoints of the
common perpendicular constructed above being close to $\ga_-\Omega^-$
and $\ga_+\Omega^+$, which have measure $0$ boundary.

\bigskip We end this section by completing the list of examples given
in \S \ref{sec:counting},  adding the following two cases.
They follow (see \cite{ParPau13b}) by applying the main Theorem
\ref{theo:mainsurvey}, the remarks following the statement of
Proposition \ref{prop:calcskinning}, and Equation
\eqref{eq:relatbowmargliouv}.

\bcoro \label{coro:bihorobitotgeod}
Let $M$ be a finite volume complete hyperbolic manifold of
dimension $n\geq 2$.

\noindent(1) If $C_-$ and $C_+$ are properly immersed finite volume
totally geodesic submanifolds of $M$ of dimensions $k_-$ and $k_+$ in
$[1,n-1]$, respectively, then, as $s\to+\infty$,
$$
\N_{C_-,C_+}(s)\sim\frac{\Vol(\SSS^{n-k_--1})\Vol(\SSS^{n-k_+-1})}
{2^{n-1}(n-1)\Vol(\SSS^{n-1})} 
\frac{\Vol(C_-)\Vol(C_+)}{\Vol(M)}\;e^{(n-1)s}\,.
$$

\noindent(2) If $\H_-$ and $\H_+$ are Margulis cusp
neighbourhoods in $M$, then, as $s\to+\infty$,
$$
\N_{\H_-,\H_+}(s)\sim\frac{2^{n-1}(n-1)\Vol(\H_-)\Vol(\H_+)}
{\Vol(\SSS^{n-1})\Vol(M)}\;e^{(n-1)s}\,.\;\;\;\Box
$$
\ecoro

In particular, if $C_-$ and $C_+$ are closed geodesics of $M$ of
lengths $\ell_-$ and $\ell_+$, respectively, then the number $\N(s)$ of
common perpendiculars (counted with multiplicity) between $C_-$ and
$C_+$ of length at most $s$ satisfies, as $s\ra+\infty$,
$$
\N(s)\sim\frac{\pi^{\frac{n}{2}-1}(\Ga(\frac{n-1}{2}))^2}
{2^{n-2}(n-1)\,\Ga(\frac{n}{2})} \,
\frac{\ell_-\;\ell_+}{\Vol(M)}\;e^{(n-1)s}\;.
$$

When $M$ is a closed hyperbolic surface (in particular $n=2$) and
$C_-=C_+$, this formula has been obtained by Martin-McKee-Wambach
\cite{MarMcKWam11} by trace formula methods. Obtaining the case
$C_-\neq C_+$, as well as error terms, seems difficult by these methods.

\section{Spectral gaps, exponential decay of correlations
 and error terms}
\label{sec:errorterm}

Let $M$ be a nonelementary complete connected Riemannian manifold with
dimension at least $2$ and pinched sectional curvature at most $-1$
having bounded derivatives. Let $\wt M\ra M$ be a universal Riemannian
cover of $M$, with covering group $\Ga$. Let $\delta$ be the critical
exponent of $\Ga$.  We assume that the Bowen-Margulis measure $m_{\rm
  BM}$ of $M$ is finite and mixing for the geodesic flow. We denote by
$\overline{m}_{\rm BM} =\frac{m_{\rm BM}}{\|m_{\rm BM}\|}$ its
normalisation to a probability measure.

In this section, we give error terms in our main counting result, when
the geodesic flow is exponentially mixing. Recall that there are two
types of exponential mixing results.

\medskip Firstly, when $M$ is locally symmetric with finite volume,
then the boundary at infinity of $\wt M$, the strong unstable,
unstable, stable, and strong stable foliations of $M$ are smooth.
Hence talking about $\C^\ell$-smooth leafwise defined functions on
$T^1M$ makes sense. We will denote by $\C_c^\ell(T^1M)$ the vector
space of $\C^\ell$-smooth functions on $T^1M$ with compact support and
by $\|\psi\|_\ell$ the Sobolev $W^{\ell,2}$-norm of any
$\psi\in\C_c^\ell(T^1M)$. Note that now the Bowen-Margulis measure of
$T^1M$ is the unique (up to a multiplicative constant) locally
homogeneous smooth measure on $T^1M$ (hence it coincides with the Liouville measure up to a
multiplicative constant which we computed in \S \ref{sec:finitevolhypman} in constant curvature).

Given $\ell\in \NN$, we will say that the geodesic flow on $T^1M$ is
{\em exponentially mixing for the Sobolev regularity $\ell$} (or that
it has {\em exponential decay of $\ell$-Sobolev correlations}) if
there exist $c,\kappa>0$ such that for all $\phi,\psi\in
\C_c^\ell(T^1M)$ and $t\in\RR$, we have
$$
\Big|\int_{T^1M} \phi\circ g^{-t}\;\psi\;d\overline{m}_{\rm BM}-\int_{T^1M}
\phi\;d\overline{m}_{\rm BM}\int_{T^1M} \psi\;d\overline{m}_{\rm BM}\;\Big|\leq
c\,e^{-\kappa |t|}\;\|\psi\|_\ell\;\|\phi\|_\ell\;.
$$
When $\Ga$ is a torsion free arithmetic lattice in the isometry group
of $\wt M$, this property, for some $\ell\in\NN$, follows from
\cite[Theo.~2.4.5]{KleMar96}, with the help of
\cite[Theo.~3.1]{Clozel03} to check its spectral gap property, and of
\cite[Lem.~3.1]{KleMar99} to deal with finite cover problems.

Secondly, when $\wt M$ is assumed to be as in the beginning of this
section, then the boundary at infinity, the strong unstable, unstable,
stable, and strong stable foliations are only H\"older smooth (as
explained in \S \ref{sec:discrisomgroups}), hence the appropriate
regularity on functions on $\wt M$ is the H\"older one.  For every
$\alpha\in\;]0,1[$, we denote by $\operatorname{C}_{\rm c} ^\alpha
(X)$ the space of $\alpha$-H\"older-continuous real-valued functions
with compact support on a metric space $(X,d)$, endowed with the
H\"older norm
$$
\|f\|_\alpha=
\|f\|_\infty+\sup_{x,\,y\in X,\;x\neq y}\frac{|f(x)-f(y)|}{d(x,y)^\alpha}\,.
$$

Given $\alpha\in\;]0,1[$, we will say that the geodesic flow
on $T^1M$ is {\em exponentially mixing for the H\"older regularity
  $\alpha$} (or that it has {\em exponential decay of
  $\alpha$-H\"older correlations}) if there exist $\kappa,c >0$ such
that for all $\phi,\psi\in \operatorname{C}_{\rm c}^\alpha(T^1M)$ and
$t\in\RR$, we have
$$
\Big|\int_{T^1M}\phi\circ\flow{-t}\;\psi\;d\overline{m}_{\rm BM}-
\int_{T^1M}\phi\; d\overline{m}_{\rm BM} 
\int_{T^1M}\psi\;d\overline{m}_{\rm BM}\;\Big|
\le c\;e^{-\kappa|t|}\;\|\phi\|_\alpha\;\|\psi\|_\alpha\,.
$$
This holds for compact manifolds $M$ when $M$ is two-dimensional by
\cite{Dolgopyat98}, when $M$ is $1/9$-pinched by
\cite[Coro.~2.7]{GiuLivPol12}, when $m_{\rm BM}$ is the Liouville
measure by \cite{Liverani04}, and when $M$ is locally symmetric by
\cite[Coro.~1.5]{Stoyanov11}.

Using smoothening processes of the functions $\phi_\eta^\pm$
introduced in the sketch of proof of \S \ref{sec:mainstatement},
we obtain the following error terms in our main counting result Theorem
\ref{theo:mainsurvey}.

\btheo[Parkkonen-Paulin \cite{ParPau13b}]\label{theo:mainsurverrterm}
Let $C_-$ and $C_+$ be two properly immersed clo\-sed convex subsets
of $M$. Assume that their skinning measures $\sigma_{C_-}$ and
$\sigma_{C_+}$ are finite and nonzero. Assume either that $M$ is
compact and the geodesic flow is exponentially mixing for the H\"older
regularity, or that $M$ is locally symmetric with finite volume and the
geodesic flow is exponentially mixing for the Sobolev regularity.
Then there is some $\kappa>0$ such that, as $s\to+\infty$,
$$
\N_{C_-,\,C_+}(s)= \frac{\|\sigma_{C_-}\|\,\|\sigma_{C_+}\|}
{\delta\,\|m_{\rm BM}\|}\;e^{\delta\, s}\,
\big(1+\operatorname{O}(e^{-\kappa\, s})\big)\;.
$$
\etheo

This error term is also valid for the effective counting Theorem
\ref{theo:mainsurveyeffective}. This result gives in particular the
exponential control in the error terms in the list of examples given
in \S \ref{sec:counting}, as well as in Corollary
\ref{coro:bihorobitotgeod}.

As an application of Theorem \ref{theo:mainsurverrterm},  using
Humbert's formula and the area of the fundamental domain 
of $\OOO_K$ in $\RR^2$ (see for example \cite[p.~318]{ElsGruMen98}),
we get a version of Cosentino's asymptotic estimate
\eqref{eq:cosentino} on the number of common perpendiculars from the
Margulis cusp neighbourhood corresponding to the horoball of points
with vertical coordinates at most $1$ to itself in $\PSL(\OOO_K)\bs
\htr$, valid for all discriminants:
$$
\N(s)=\frac{\pi\,|\OOO_K^\times|^2}{4\sqrt{|D_K|}\,\zeta_K(2)}
\;e^{2s}\big(1+\operatorname{O}(e^{-\kappa s})\big),
$$
when $s\to+\infty$. We refer to \cite{ParPauTou} for further
generalisations.

\section{Gibbs measures and counting arcs with weights}
\label{sec:gibbs}

Let $M$ be a nonelementary complete connected Riemannian manifold with
dimension at least $2$ and pinched sectional curvature at most
$-1$. Let $\wt M\ra M$ be a universal Riemannian cover of $M$, with
covering group $\Ga$, and $x_0\in\wt M$.

When counting geodesic arcs, it is sometimes useful to give them a
higher weight if they are passing through a given region of $M$, and
even more precisely, through a given region in position and
direction. The trick is to introduce a {\em potential}, that is a
H\"older-continuous map $F:T^1M\ra\RR$. To shorten the exposition, we
will assume in this survey that $F$ is bounded and {\em reversible},
that is, that $F\circ\iota=F$. These assumptions are not necessary, up
to the appropriate modifications, see \cite{ParPau13b}. Given a
piecewise smooth path $c:[a,b]\ra M$, one defines its {\em weighted
  length} for the potential $F$ as
$$
\int_c F=\int_a^{b}F\circ\dot c\,(t)\;dt\;.
$$
We are now going to adapt the material of \S \ref{sec:skinning}
to the weighted case, see for instance \cite{Ledrappier95,
  Hamenstadt97, Coudene03, Schapira04a, Mohsen07, PauPolSha,
  ParPau13b} with an emphasis on the last two ones for more
information.

Let $\wt F=F\circ Tp:T^1\wt M\ra\RR$ be the lift of $F$ by the
universal cover $p:\wt M\ra M$. For every $x,y$ in $\wt M$, if
$c:[0,d(x,y)]\ra\wt M$ is the geodesic path from $x$ to $y$, let
$\int_x^y\wt F=\int_0^{d(x,y)}\wt F\circ \dot c\,(t) \;dt$. The {\em
  critical exponent} of the potential $F$ is
$$
\delta_{F}=\lim_{n\to+\infty}\frac
1n\ln\sum_{\ga\in\Ga, \,d(x,\ga y)\le n}e^{\int_x^{\ga y}\wt F}\;,
$$
see \cite[Theo.~4.2]{PauPolSha} for the existence and finiteness of
the above limit and its independence on $x,y\in \wt M$. Replacing the
previous critical exponent $\delta$ (to which it is equal if $F= 0$),
the critical exponent $\delta_F$ of the potential $F$ will give the
exponential growth rate in the counting of weighted common
perpendiculars.

Similarly, the Busemann cocycle $\beta_\xi(x,y)$ needs to be
replaced. The (normalized) {\em Gibbs cocycle} associated to the
potential $F$ is the function $C=C^F:\partial_\infty\wt M\times\wt
M\times\wt M\to\RR$ defined by
$$
(\xi,x,y)\mapsto C_\xi(x,y)=\lim_{t\to+\infty}
\int_y^{\xi_t}(\wt F-\delta_F)- \int_x^{\xi_t}(\wt F-\delta_F)\,,
$$
where $t\mapsto\xi_t$ is any geodesic ray with endpoint
$\xi\in\partial_\infty\wt M$.  The Gibbs cocycle is well defined by
the H\"older-continuity of $F$. It satisfies obvious equivariance and
cocycle properties: For all $x,y,z\in\wt M$, and for every isometry
$\ga$ of $\wt M$, we have
\begin{equation}\label{eq:cocycleGibbs}
 C_{\ga \xi}(\ga x,\ga y)= C_\xi(x,y)\; \textrm{ and }\;
 C_\xi(x,z)+ C_\xi(z,y)= C_\xi(x,y)\;.
\end{equation}

Similarly, the Bowen-Margulis measure needs to be replaced. A family
$(\mu_x)_{x\in\wt M}$ of finite Borel measures on $\partial_\infty\wt
M$, whose support is the limit set $\Lambda\Ga$ of $\Ga$, is a {\em
  Patterson density for the potential $F$} (of dimension $\delta_F$)
if
$$
\ga_*\mu_x=\mu_{\ga x}
$$
for all $\ga\in\Ga$ and $x\in \wt M$, and if the following
Radon-Nikodym derivative exists for all $x,y\in\wt M$ and satisfies,
for all $\xi\in\partial_\infty\wt M$,
$$
\frac{d\mu_x}{d\mu_y}(\xi)=e^{-C_\xi(x,\,y)}\;.
$$
Let $(\mu_x)_{x\in\wt M}$ be such a Patterson density. The {\em Gibbs
  measure} on $T^1\wt M$ for the potential $F$ is the measure $\wt
m_{F}$ on $T^1\wt M$ given by
$$
d\wt m_F(v)=
e^{C_{v_-}(x_0,\,\pi(v))+C_{v_+}(x_0,\,\pi(v))}\,
d\mu_{x_0}(v_-)\,d\mu_{x_0}(v_+)\,dt\;,
$$ 
using Hopf's parametrisation. The Gibbs measure $\wt m_F$ is
independent of the base point ${x_0}\in\wt M$ used in its definition,
and it is invariant under the actions of the group $\Ga$ and the
geodesic flow. Thus, it defines a measure $m_{F}$ on $T^1M$ which is
invariant under the geodesic flow, called the {\em Gibbs measure} on
$T^1M$ for the potential $F$. When the Gibbs measure $m_{F}$ is
finite, there exists a unique (up to a multiplicative constant)
Patterson density for the potential $F$; the
probability measure $\frac{m_{F}}{\|m_{F}\|}$ is uniquely
defined; it is the unique probability measure of maximal pressure for
the geodesic flow and the potential $F$; see \cite{PauPolSha} for
proofs of these claims. When finite, the Gibbs measure on $T^1M$ is
mixing if the geodesic flow is topologically mixing, see
\cite{Babillot02b}.

Let $D$ be a properly immersed closed convex subset of $M$, and let
$\wt D$ be a proper nonempty closed convex subset of $\wt M$, whose
$\Ga$-orbit is locally finite and whose image in $M$ is $D$. We also
need to adapt the skinning measures to the presence of the potential
$F$. The {\em skinning measure} of $\wt D$ for the potential $F$ is
the measure $\wt\sigma^F_{\wt D}$ on $\partial^1_+ \wt D$, defined,
using the homeomorphism $v\mapsto v_+$ from $\partial^1_+ \wt D$ to
$\partial_\infty \wt M-\partial_\infty \wt D$, by
$$
d\wt\sigma^F_{\wt D}(v) = 
e^{C_{v_+}(x_0,\,P_{\wt D}(v_+))}\;d\mu_{x_0}(v_{+})\;.
$$
It is independent of the base point $x_0$, and satisfies
$\ga_*(\wt\sigma^F_{\wt D})=\wt\sigma^F_{\ga \wt D}$ for every
$\ga\in\Ga$. Let $\Ga_\wt D$ be the stabiliser in $\Ga$ of $\wt
D$. The $\Ga$-invariant locally finite Borel positive measure
$\sum_{\ga\in\Ga/\Ga_{\wt D}} \ga_*\wt\sigma^F_{\wt D}$ defines,
through the covering $T^1\wt M\ra T^1M$, a locally finite measure
$\sigma^F_D$, called the {\em skinning measure} of $D$ for the
potential $F$. See \cite{ParPau13b} for further information on the
skinning measures with potential. 

Let $D_-,D_+$ be two properly immersed closed convex subsets of
$M$. For every $s\geq 0$, recall that $\Perp_{D_-,D_+}(s)$ is the set
of the common perpendiculars from $D_-$ to $D_+$ having lengths at most
$s$. The {\em weighted counting function} of common perpendiculars
between $D_-$ and $D_+$ (counted with multiplicities) for the
potential $F$ is
$$
\N_{D_-,D_+,F}(s)=\sum_{c\in \Perp_{D_-,D_+}(s)} m(c) \;e^{\int_cF}\;.
$$
Using the same scheme of proof as explained in \S
\ref{sec:mainstatement}, we have the following asymptotic result.

\btheo[Parkkonen-Paulin \cite{ParPau13b}]\label{theo:mainsurveyGibbs}
Let $M$ be a nonelementary complete connected Riemannian manifold with
pinched sectional curvature $-b^2\le K\le -1$, and let $F:T^1M\ra \RR$
be a (bounded,  reversible) H\"older-continuous map. Let
$\delta_F$ be the critical exponent of the potential $F$.  Assume that
the Gibbs measure $m_F$ is finite and mixing for the geodesic flow.
Let $D_-$ and $D_+$ be two properly immersed closed convex subsets of
$M$. Assume that $\sigma^F_{D_-}$ and $\sigma^F_{D_+}$ are
finite and nonzero. Then, as $s\ra +\infty$,
$$
\N_{D_-,D_+,F}(s)\sim \frac{\|\sigma^F_{D_-}\|\,\|\sigma^{F}_{D_+}\|} 
{\delta_F\,\|m_F\|}\;e^{\delta_F s}\,.
$$
\etheo

We have error terms in the presence of exponential decay of
correlations, and the endpoints of the common perpendiculars are
evenly 
distributed, that is, we may restrict to counting the common
perpendiculars with endpoints in measurable subsets $\Omega_-$ and
$\Omega_+$, with finite nonzero skinning measures for the potential
$F$ and negligible boundary, of $\partial^1_+D_-$ and
$\partial^1_+D_+$, respectively. As in the two previous sections, we refer to 
\cite{ParPau13b} for precise statements and proofs.

 {\small
  \bibliography{../biblio} }
%{\small \bibliography{../viitteet} }

\bigskip
{\small\noindent \begin{tabular}{l} 
Department of Mathematics and Statistics, P.O. Box 35\\ 
40014 University of Jyv\"askyl\"a, FINLAND.\\
{\em e-mail: jouni.t.parkkonen@jyu.fi}
\end{tabular}
\medskip

\noindent \begin{tabular}{l}
D\'epartement de math\'ematique, UMR 8628 CNRS, B\^at.~425\\
Universit\'e Paris-Sud,
91405 ORSAY Cedex, FRANCE\\
{\em e-mail: frederic.paulin@math.u-psud.fr}
\end{tabular}}

\end{document}